\documentclass[a4paper,final]{article}
\usepackage{amsfonts,latexsym,url,amsmath,enumitem,amsmath}

\usepackage[color,notcite]{showkeys}
\definecolor{labelkey}{rgb}{0.6,0,1}

\newtheorem{theorem}{Theorem}[section]
\newtheorem{prop}[theorem]{Proposition}
\newtheorem{lem}[theorem]{Lemma}
\newtheorem{cor}[theorem]{Corollary}
\newtheorem{defi}[theorem]{Definition}
\newtheorem{remark}[theorem]{Remark}

\def\bt{\begin{theorem}}
\def\et{\end{theorem}}
\def\bp{\begin{prop}}
\def\ep{\end{prop}}
\def\bl{\begin{lem}}
\def\el{\end{lem}}
\def\bc{\begin{cor}}
\def\ec{\end{cor}}
\def\bd{\begin{defi}}
\def\ed{\end{defi}}
\def\br{\begin{remark}}
\def\er{\end{remark}}

\def\findem{\rule{0.2cm}{0.2cm}}
\def\bpr#1{\textbf{Proof#1}}
\def\epr{\findem}

\def\bhyp#1{\begin{equation}\label{#1}\begin{array}{c}} 
\def\ehyp{\end{array}\end{equation}}

\def\bpro#1{\begin{equation}\label{#1}\left\{\begin{array}{l@{\qquad}l}}
 
\def\epro{\end{array}\right.\end{equation}}

\newcounter{cst}
\def \ctel#1{C_{\refstepcounter{cst}\label{#1}\thecst}}
\def \cter#1{C_{\ref{#1}}}

\def\listeeq#1{\mbox{\rm [#1]}}

\def\Annexe{
        \appendix
        
        \setcounter{equation}{0}
        \setcounter{theorem}{0}
        }

\newtheorem{theoreme-Annexe}{Th\'eor\`eme}[section]
\newtheorem{prop-Annexe}[theoreme-Annexe]{Proposition}
\newtheorem{cor-Annexe}[theoreme-Annexe]{Corollaire}
\newtheorem{lem-Annexe}[theoreme-Annexe]{Lemme}
\newtheorem{defi-Annexe}[theoreme-Annexe]{D\'efinition}
\newtheorem{remarque-Annexe}[theoreme-Annexe]{Remarque}

\def\R{\mathbf{R}}

\def\n{\mathbf{n}}

\def\dx{\,{\rm d}x}

\def\<{\langle}
\def\>{\rangle}

\def\dsp{\displaystyle} 
\def\div{{\rm div}}
\def\refe#1{\eqref{#1}}

\def\divmfd{\mathcal{DIV}^h}
\def\gradmfd{\mathcal{G}^h}

\def\Xm{H_{\mathcal M}}
\def\Xmb{\widetilde{H}_{\mathcal M}}
\def\Xmbzero{\widetilde{H}_{\mathcal M,0}}
\def\XmbzeroB{\widetilde{H}_{\mathcal M,0}^{\mathcal B}}
\def\Xmbloc#1{\widetilde{H}_{#1}}
\def\ploc#1{\widetilde{p}_{#1}}
\def\qloc#1{\widetilde{q}_{#1}}

\def\v{\mathbf{v}}
\def\xs{\bar x_{\sigma}}

\begin{document}

\begin{center} {\Large\textbf{A unified approach to Mimetic Finite Difference, Hybrid Finite Volume and
Mixed Finite Volume methods}}\footnote{This work was supported by GDR MOMAS CNRS/PACEN}

\vspace*{0.7cm}

 J\'er\^ome Droniou \footnote{Universit\'e Montpellier 2, Institut de Math\'ematiques et de Mod\'elisation de
Montpellier, CC 051,
Place Eug\`ene Bataillon, 34095 Montpellier cedex 5, France. email: \texttt{droniou@math.univ-montp2.fr}},
Robert Eymard \footnote{Universit\'e Paris-Est, Laboratoire d'Analyse et de Math\'ematiques Appliqu\'ees,
UMR 8050, 5 boulevard Descartes, Champs-sur-Marne, 77454 Marne-la-Vall\'ee Cedex 2,
France. email:\texttt{Robert.Eymard@univ-mlv.fr}},
Thierry Gallou\"et\footnote{L.A.T.P., UMR 6632, Universit\'e de Provence, \texttt{gallouet@cmi.univ-mrs.fr}},
Rapha\`ele Herbin\footnote{L.A.T.P., UMR 6632, Universit\'e de Provence, \texttt{herbin@cmi.univ-mrs.fr}}.

\today

\end{center}

\textbf{Abstract}
We investigate the connections between several recent methods for the discretization of  ani\-so\-tropic heterogeneous diffusion operators on general grids. 
We prove that the Mimetic Finite Difference scheme, the Hybrid Finite Volume scheme and the Mixed Finite Volume scheme are in fact identical up to some slight generalizations. 
As a consequence,  some of the mathematical results obtained for each of the method (such as convergence properties or error estimates) may be extended to the unified common framework.
We then focus on the relationships between this unified method and nonconforming Finite Element schemes or Mixed Finite Element schemes, obtaining as a by-product an explicit lifting operator close to the ones used in some theoretical studies of the Mimetic Finite Difference scheme.
We also show that for isotropic operators, on particular meshes such as triangular meshes with acute angles, the unified method boils down to the well-known efficient two-point flux Finite Volume scheme.

\section{Introduction}\label{sec-intro}

A benchmark was organized at the last FVCA 5 conference \cite{her-08-ben} in June 2008 to test the recently developed 
methods for the numerical solution of heterogeneous anisotropic problems.
In this benchmark and in this paper, we consider the Poisson equation with homogeneous boundary condition
\begin{subequations}
\begin{align}
-\div(\Lambda \nabla p) = f   &\quad\mbox{ in $\Omega$}, \label{baseq}\\
p=0 &\quad\mbox{ on $\partial\Omega$},                   \label{basecl}
\end{align}
\label{base}\end{subequations}
where $\Omega$ is a bounded open subset of $\R^d$ ($d\ge 1$), $\Lambda:\Omega\to M_d(\R)$ 
is bounded measurable symmetric and uniformly elliptic (i.e. there exists $\zeta>0$ such that, 
for a.e. $x\in\Omega$ and all $\xi\in\R^d$, $\Lambda(x)\xi\cdot\xi\ge \zeta|\xi|^2$) and $f\in L^2(\Omega)$.

The results of this benchmark, in particular those of  \cite{cha-08-ben,eym-08-ben,lip-08-ben,man-08-ben} seem to demonstrate that the behavior of three of the submitted methods, namely the Hybrid Finite Volume method \cite{eym-07-new,eym-08-dis2,eym-08-dis}, the  Mimetic Finite Difference method \cite{bre-05-con,bre-05-fam}, and the Mixed Finite Volume method \cite{dro-06-mix,dro-08-stu} are quite similar in a number of cases (to keep notations light while retaining the legibility, in the following we call ``Hybrid'', ``Mimetic'' and ``Mixed''  these  respective methods).

\medskip

A straightforward common point between these methods is that they are written using a  general partition of $\Omega$ into 
polygonal open subsets and that they introduce unknowns which approximate the solution $p$ and the fluxes of its gradient 
on the edges of the partition. 
However, a comparison of the methods is still lacking, probably because their mathematical analysis relies on different 
tool boxes. 
The mathematical analysis of the Mimetic method \cite{bre-05-con} is based on an error estimate technique (in the spirit 
of the mixed finite element  methods). 
For the Hybrid and the Mixed methods \cite{eym-08-dis,dro-06-mix}, the convergence proofs rely on discrete functional 
analysis tools.
The aim of this paper is to point out the common points of these three methods. 
To this purpose, we first gather, in Table \ref{tabnot} of  Section \ref{sec-methods}, some definitions and notations 
associated with each method, and we present the three methods as they are introduced in the literature;  we also present a
generalized or modified form of each of the methods. The three resulting methods are then shown to be 
identical (Section \ref{algebraic}) and to inherit some of the mathematical properties of the initial methods
(Section \ref{sec-convergence}).
Particular cases are  then explored in Section \ref{sec-links}, which show the thorough relations between this unified 
method and a few classical methods (nonconforming Finite Elements,  Mixed Finite Elements and two-point flux Finite 
Volumes) and, as a by-product, provide a flux lifting operator for a particular choice of the stabilization parameters 
(but any kind of grid).
Finally, a few technical results are provided in an appendix.

\medskip

In this article, we only consider the case of a linear single equation \eqref{base}, but the three methods we study have 
also been used on more complex problems, such as the incompressible Navier-Stokes equations
\cite{che-09-col,dro-08-stu}, fully non-linear equations of the $p$-Laplacian type \cite{dro-06-plap,eym-09-cel}, 
non-linear coupled problems \cite{cha-07-por}, etc.
However, the main ideas to apply these methods in such more complex situations stem from the study of their properties on 
the linear diffusion equation.
The unifying framework which is proposed here has a larger field of application than \eqref{base}; it facilitates the 
transfer of  the ideas and techniques used for one method  to another and can also  give some new  leads for each method.

\section{The methods}\label{sec-methods}

We first provide the definitions and notations associated with each method (in order for the readers who are familiar with one or the other theory to easily follow the subsequent analysis, we shall freely use one or the other notation (this of course yields some redundant notations).
\setcounter{equation}{0}
\begin{table}[h!]
\begin{center}
\begin{tabular}{|l|l|l|}
\hline
\rule{0em}{1.1em}&\emph{Mimetic notation}&\emph{Finite Volume notation}\\[0.1em]
\hline
\rule{0em}{1.1em}Partition of $\Omega$ in polygonal sets& $\Omega_h$ & $\mathcal M$\\[0.1em]
\hline
\rule{0em}{1.1em}Elements of the partition (grid cells) & $E$ & $K$ (``control volume'')\\[0.1em]
\hline
\rule{0em}{1.1em}Set of edges/faces of a grid cell &
\begin{minipage}{8em}\rule{0em}{1.1em}$\partial E$,
or numbered from $1$ to $k_E$\end{minipage} & $\mathcal E_K$\\[0.1em]
\hline
\rule{0em}{1.1em}Edges/faces of grid cell & $e$ & $\sigma$\\[0.1em]
\hline
\begin{minipage}{22em}\rule{0em}{1.1em}Space of discrete $p$ unknowns
(piecewise approximations of $p$ on the partition)\end{minipage} & $Q^h$ & $\Xm$\\[0.1em]
\hline
\begin{minipage}{22em}\rule{0em}{1.1em}Approximation of the solution $p$ on
the grid cell $E=K$\end{minipage}
& $p_E$ & $p_K$\\[0.1em]
\hline
\begin{minipage}{22em}\rule{0em}{1.1em}Discrete flux: approximation of $\frac{1}{|e|}\int_{e}
-\Lambda \nabla p\cdot\n_{E}^e$
(with $e=\sigma$ an edge of $E=K$ and $\n_{E}^e=\n_{K,\sigma}$ the unit normal to $e$ outward $E$)\end{minipage}
& $F^e_E$ & $F_{K,\sigma}$\\
\hline
\begin{minipage}{22em}\rule{0em}{1.1em}Space of discrete fluxes
(approximations of $\frac{1}{|e|}\int_{e/\sigma} -\Lambda \nabla p\cdot\n$)\end{minipage}
& $X^h$ & $\mathcal F$\\[0.1em]
\hline
\end{tabular}
\caption{Usual notations and definitions in Mimetic and Finite Volume frameworks.\label{tabnot}}
\end{center}
\end{table}

\br In the usual Finite Volume literature, the quantity $F_{K,\sigma}$ is usually rather an approximation of $\int_{\sigma}  -\Lambda \nabla p\cdot\n_{K,\sigma}$  than $\frac{1}{|\sigma|}\int_{\sigma} -\Lambda \nabla 
p\cdot\n_{K,\sigma}$;
we choose here the latter normalization in order to simplify the comparison.
\er

In all three methods, a natural condition is imposed on the gradient fluxes
(called conservativity in Finite Volume methods and continuity condition in Mimetic methods):
for any interior edge  $\sigma$ (or $e$) between two polygons $K$ and $L$ (or $E_1$ and
$E_2$),
\begin{equation}
F_{K,\sigma}+F_{L,\sigma}=0,\quad\mbox{ (or } F_{E_1}^e=-F_{E_2}^e).
\label{conservativite}\end{equation}
This condition is included in the definition of the discrete flux space $\mathcal F$ (or  $X^h$).

\medskip

The Mimetic, Hybrid and Mixed methods for \refe{base} all consist in seeking $p\in\Xm$ (or $ Q^h$) and $F\in \mathcal F$ (or  $X^h$), which approximate respectively the solution $p$ and its gradient fluxes,
writing equations on these unknowns which discretize the continuous equation \refe{base}. 
Each method is in fact a family of schemes rather than a unique scheme: indeed, there exists some freedom on the choice of some of the  parameters  of the scheme (for instance in the stabilization terms which ensure the coercivity of the methods).

\subsection{The Mimetic method}\label{sec-mfd}

The standard Mimetic method first consists in defining, from the Stokes formula, 
a discrete divergence operator on the space of the discrete fluxes:
for $G\in X^h$, $\divmfd G\in Q^h$ is defined by
\begin{equation}
(\divmfd G)_E=\frac{1}{|E|}\sum_{i=1}^{k_E} |e_i| G_E^{e_i}.
\label{defdivmfd}
\end{equation}
The space $Q^h$ of piecewise constant functions is endowed with the usual $L^2$ inner product
\[
[p,q]_{Q^h}=\sum_{E\in\Omega_h}|E|p_Eq_E
\]
and a local inner product is  defined on the space of fluxes unknowns of each element $E$: 
\begin{equation}
[F_E,G_E]_E= F_E^t \mathbb{M}_{E} G_E = \sum_{s,r=1}^{k_E}\mathbb{M}_{E,s,r} F_E^{e_s} G_E^{e_r},
\label{localscal}\end{equation}
where $ \mathbb{M}_{E} $ is a symmetric definite positive matrix of order $k_E$.
Each local inner product is assumed to satisfy the following discrete Stokes formula (called Condition  ({\bf S2}) in \cite{bre-05-con,bre-05-fam}):
\begin{equation}
\dsp\forall E\in\Omega_h\,,\;\forall q\mbox{ affine function, }\forall G\in X^h\,:\;
[(\Lambda\nabla q)^I,G]_E+\int_E q \ (\divmfd G)_E\,dV=\sum_{i=1}^{k_E}G_E^{e_i}\int_{e_i}q\,d\Sigma
\label{S2}\end{equation}
where $((\Lambda\nabla q)^I)_E^e=\frac{1}{|e|}\int_e \Lambda_E\nabla q\cdot\n_{E}^e\,d\Sigma$ and  $\Lambda_E$ is the value, assumed to be constant, of $\Lambda$ on $E$ (if $\Lambda$ is not constant on $E$, one can take $\Lambda_E$ equal to the mean value of $\Lambda$ on $E$).
\begin{remark}
Note that Condition  ({\bf S1}) of \cite{bre-05-con,bre-05-fam} is needed in the convergence study of the method, but not in its definition; therefore it is only recalled in Section \ref{sec-convergence}, see \refe{S1}.
\end{remark}
The local inner products \refe{localscal} allow us to construct a complete inner product $[F,G]_{X^h}=\sum_{E\in\Omega_h}[F,G]_E$, which in turn defines a discrete flux operator $\gradmfd:Q^h\to X^h$ as the adjoint operator of $\divmfd$: for all $F\in X^h$ and $p\in Q^h$,
\[
[F,\gradmfd p]_{X^h}=[p,\divmfd F]_{Q^h}
\]
(notice that this definition of $\gradmfd p$ takes into account the homogeneous boundary condition \refe{basecl}).
Using these definitions and notations, the standard Mimetic method then reads:
find $(p,F)\in Q^h\times X^h$ such that
\begin{equation}
\divmfd F=f_h\,,\quad F=\gradmfd p
\label{mfdscheme}\end{equation}
where $f_h$ is the $L^2$ projection of $f$ on $Q^h$, or in the equivalent weak form:
find $(p,F)\in Q^h\times X^h$ such that
\begin{eqnarray}
&&\forall G\in X^h\,:\;[F,G]_{X^h}-[p,\divmfd G]_{Q^h}=0\,,\label{defmfd1}\\
&&\forall q\in Q^h\,:\;[\divmfd F,q]_{Q^h}=[f_h,q]_{Q^h}.\label{defmfd2}
\end{eqnarray}
The precise definition of the Mimetic method requires to choose the local matrices $\mathbb{M}_E$ defining the local inner products $[\cdot,\cdot]_E$.
It can be shown  \cite{bre-05-fam} (see also Lemma \ref{lem-equivGS2} in the appendix) that this matrix defines an inner product satisfying \refe{S2} if and only if it can be written $ \mathbb{M}_E   \mathbb{N}_E = \bar{\mathbb{R}}_E $ or equivalently
\begin{equation}
\mathbb{M}_E=\frac{1}{|E|}\bar{\mathbb{R}}_E\Lambda_E^{-1}\bar{\mathbb{R}}_E^T
+\mathbb{C}_E\mathbb{U}_E\mathbb{C}_E^T
\label{choixmatmfd-initial}\end{equation}
where
\begin{equation}
\begin{array}{c}
\bar{\mathbb{R}}_E \mbox{ is the } k_E\times d \mbox{ matrix with rows }(|e_i|(\bar x_{e_i}-\bar x_E)^T)_{i=1,k_E}, \\
\mbox{with } \bar x_e \mbox{ the center of gravity of the edge } e \mbox{ and } \bar x_E \mbox{ the center of gravity of the cell } E,
\end{array}
\label{mf1}
\end{equation}

\begin{equation} \left.\begin{array}{c} \mathbb{C}_E\mbox{ is a $k_E\times (k_E-d)$ matrix such that
${\rm Im}(\mathbb{C}_E)=({\rm Im}(\mathbb{N}_E))^\bot$,}\\
\mbox{ with $\mathbb{N}_E$ the } k_E\times d \mbox{ matrix with columns }\\
(\mathbb{N}_E)_j=\left(\begin{array}{c} (\Lambda_E)_j\cdot \n_E^{e_1}\\
\vdots\\ 
(\Lambda_E)_j\cdot \n_E^{e_{k_E}}\end{array}\right)\quad\mbox{for $j=1,\ldots,d$}, \\ \mbox{$(\Lambda_E)_j$
being the $j$-th column of $\Lambda_E$}
\end{array}\right\}
\label{mf2}
\end{equation}
\begin{equation}
\mathbb{U}_E\mbox{ is a } (k_E-d)\times (k_E-d)\mbox{ symmetric definite positive matrix}.
\label{mf3}
\end{equation}

\medskip

Here, we consider a slightly more general version of the Mimetic method, replacing $\bar x_E$  by a point $x_E$ which may be chosen different from the center of gravity of $E$. 
We therefore take
\begin{equation}
\mathbb{M}_E=\frac{1}{|E|}\mathbb{R}_E\Lambda_E^{-1}\mathbb{R}_E^T
+\mathbb{C}_E\mathbb{U}_E\mathbb{C}_E^T
\label{choixmatmfd}\end{equation}
where
\begin{equation}\label{defnewR}
\begin{array}{c}
\mbox{$\mathbb{R}_E$ is the $k_E\times d$ matrix with rows
$(|e_i|(\bar x_{e_i}-x_E)^T)_{i=1,k_E}$,}\\
\mbox{with $\bar x_e$ the center of gravity of the edge $e$ and
$x_E$ any point in the cell $E$.}
\end{array}
\end{equation}
The other matrices $\mathbb{C}_E$ and $\mathbb{U}_E$ remain given by \eqref{mf2} and \eqref{mf3}.
The choice \listeeq{\eqref{choixmatmfd},\eqref{defnewR}} of $\mathbb{M}_E$ no longer gives, in general, an inner product 
$[\cdot,\cdot]_E$ which satisfies \refe{S2}, but it yields a generalization of this assumption; indeed,
choosing a weight function $w_E:E\to \R$ such that
\begin{equation}
\int_E w_E(x)\dx = |E|\quad\mbox{ and }\quad \int_E xw_E(x)\dx=|E|x_E,
\label{hypowE}\end{equation}
we prove in the appendix (Section \ref{appen-genmim}) that the matrix $\mathbb{M}_E$ can be written \refe{choixmatmfd} with $\mathbb{R}_E$ defined by \refe{defnewR} if and only if the corresponding inner product $[\cdot,\cdot]_E$ satisfies
\begin{equation}
\begin{array}{ll}\dsp\forall E\in\Omega_h\,,\;\forall q\mbox{ affine function, }\forall G\in X^h\,:&
\displaystyle [(\Lambda\nabla q)^I,G]_E+\int_E q (\divmfd G)_E\  w_E\,dV\\
&\displaystyle =\sum_{i=1}^{k_E}G_E^{e_i}\int_{e_i}q\,d\Sigma.
\end{array}
\label{GS2}\end{equation}

\begin{defi}[Generalized Mimetic method] The Generalized Mimetic scheme for \refe{base} reads:
\begin{center}
Find $(p,F)\in Q^h\times X^h$ which satisfies
\listeeq{\refe{defdivmfd},\refe{localscal},\refe{defmfd1},\refe{defmfd2},\refe{choixmatmfd},\refe{defnewR}}.
\end{center}
Its parameters are the family of points $(x_E)_{E\in \Omega_h}$ (which are freely chosen inside each grid cell)
and the family of stabilization matrices $(\mathbb{C}_E,\mathbb{U}_E)_{E\in\Omega_h}$
satisfying \eqref{mf2} and \eqref{mf3}.
\label{def-genmimetic}\end{defi}

\subsection{The Hybrid method}\label{sec-hfv}

The standard Hybrid method is best defined using additional unknowns $p_\sigma$ playing the role of approximations of $p$ on the edges of the discretization of $\Omega$; if $\mathcal E$ is the set of such edges, we let $\Xmb$ be the extension of $\Xm$ consisting of vectors $\widetilde{p}=((p_K)_{K\in\mathcal M}, (p_\sigma)_{\sigma\in\mathcal E})$. 
It will also be useful to consider the space $\Xmbloc{K}$ of the restrictions $\ploc{K}=(p_K,(p_\sigma)_{\sigma\in\mathcal E_K})$ to a control volume $K$ and its edge of the elements $\widetilde{p}\in\Xmb$.
A discrete gradient inside $K$ is defined for $\ploc{K}\in \Xmbloc{K}$ by
\begin{equation}
\nabla_K \ploc{K}=\frac{1}{|K|}\sum_{\sigma\in\mathcal E_K}|\sigma|(p_\sigma-p_K)\n_{K,\sigma}.
\label{defgradhfv}\end{equation}
The definition of this discrete  gradient stems from the  fact that, for any vector $\mathbf{e}$, any control volume $K$ and any $x_K\in \R^d$, we have
\begin{equation}
|K| \mathbf{e}=\sum_{\sigma\in\mathcal E_K} |\sigma| \mathbf{e}\cdot(\xs -x_K)\n_{K,\sigma}
\label{formmag2}
\end{equation}
where $\xs$ is the center of gravity of $\sigma$ and $x_K$ is any point of $K$.
Hence the discrete gradient is  consistent in the sense that if $p_\sigma = \psi(\overline x_\sigma)$ and $p_K = \psi(x_K)$ for an affine function $\psi$ on $K$,  then  $\nabla_K \ploc{K} = \nabla \psi $ on $K$ (note that this consistency is not sufficient to ensure strong convergence, and in fact, only weak convergence of the discrete gradient will be obtained).
A stabilization needs to be added to the discrete gradient \eqref{defgradhfv} in order to build a discrete coercive bilinear form expected to approximate the bilinear form $(u,v) \mapsto \int_\Omega \Lambda \nabla u \cdot \nabla v \dx$ occurring in the weak formulation of \eqref{base}.  
Choosing a point $x_K$ for each control volume $K$ (for instance the center of gravity, but this is not mandatory), and keeping in mind that $p_K$ is expected to be an approximation of the solution $p$ of \eqref{base} at this point, a second-order consistency error term  $S_K(\ploc{K})=(S_{K,\sigma}(\ploc{K}))_{\sigma\in\mathcal E_K}$ is defined by
\begin{equation}
S_{K,\sigma}(\ploc{K})=p_\sigma-p_K-\nabla_K \ploc{K}\cdot(\xs-x_K).
\label{defSk}\end{equation}
Note that thanks to \eqref{formmag2}, 
 \begin{equation}
\sum_{\sigma\in\mathcal E_K}|\sigma| S_{K,\sigma}(\ploc{K})\n_{K,\sigma} = 0, \label{consvhf}
\end{equation}
and that $S_{K,\sigma}(\ploc{K})=0$ if $p_\sigma = \psi(\overline x_\sigma)$ and $p_K = \psi(x_K)$ for an affine function $\psi$ on $K$.

The fluxes $(F_{K,\sigma})_{\sigma\in\mathcal E_K}$ on the boundary of a control volume $K$ associated with some $\widetilde{p}\in \Xmb$ are then defined by imposing that
\begin{eqnarray}
&&\forall K\in\mathcal M\,,\;
\forall \qloc{K}=(q_K,(q_\sigma)_{\sigma\in\mathcal E_K})\in\Xmbloc{K}\,:\nonumber\\
&&\sum_{\sigma\in\mathcal E_K} |\sigma|F_{K,\sigma}(q_K-q_\sigma)
=|K|\Lambda_K \nabla_K \ploc{K}\cdot\nabla_K \qloc{K}
+\sum_{\sigma\in\mathcal E_K}\alpha_{K,\sigma} \frac{|\sigma|}{d_{K,\sigma}}
S_{K,\sigma}(\ploc{K})S_{K,\sigma}(\qloc{K})
\label{deffluxhfv1}\end{eqnarray}
where $\Lambda_K$ is the mean value of $\Lambda$ on $K$,   $d_{K,\sigma}$ is the distance between $x_K$ and the hyperplane containing $\sigma$ and $\alpha_{K,\sigma}>0$. 
Note that $F_{K,\sigma}$ is uniquely defined by \eqref{deffluxhfv1}, since this equation is equivalent to 
\begin{equation}
|\sigma| F_{K,\sigma} = |K|\Lambda_K \nabla_K \ploc{K}\cdot\nabla_K \qloc{K}
+\sum_{\sigma\in\mathcal E_K}\alpha_{K,\sigma} \frac{|\sigma|}{d_{K,\sigma}}
S_{K,\sigma}(\ploc{K})S_{K,\sigma}(\qloc{K})
\label{fluxhfv1}\end{equation}
where $\qloc{K}$ satisfies $q_K - q_\sigma = 1$ and  $q_K - q_{\sigma'} = 0$  for $\sigma' \ne \sigma.$
To take into account the boundary condition \eqref{basecl}, we denote by $\Xmbzero = \{ \widetilde{p}\in \Xmb$ such that $p_\sigma=0$ if $\sigma\subset \partial\Omega\}$ and the Hybrid method can then be written: find $\widetilde{p}\in \Xmbzero$ such that, with $F_{K,\sigma}$ defined by \refe{deffluxhfv1},
\begin{equation}
\forall \widetilde{q}\in \Xmbzero\,:\;
\sum_{K\in\mathcal M}\sum_{\sigma\in\mathcal E_K}|\sigma|F_{K,\sigma}(q_K-q_\sigma)
=\sum_{K\in\mathcal M}q_K\int_K f.
\label{defhfv}\end{equation}
In particular, taking  $q_K=0$ for all $K$, and $q_\sigma=1$ for one $\sigma$ and $0$ for the others in \eqref{defhfv}  yields   that $F$ satisfies \refe{conservativite}; once  this conservativity property is imposed by requiring that $F \in \mathcal F$,  we may eliminate the $q_\sigma$ from \refe{defhfv} and reduce the Hybrid method to: find $(\widetilde{p},F)\in \Xmbzero\times\mathcal F$ satisfying \refe{deffluxhfv1}
and
\begin{equation}
\forall K\in\mathcal M\,:\;
\sum_{\sigma\in\mathcal E_K}|\sigma|F_{K,\sigma}
=\int_K f.
\label{defhfv_bis}\end{equation}
This last equation is the flux balance, one of the key ingredients of the finite volume methods.

\medskip

Let us now introduce a generalization of the Hybrid method  for the comparison with the other methods.
As previously mentioned, the stabilization term $S_K$ is added in \refe{deffluxhfv1} in order to obtain a coercive scheme
(the sole discrete gradient $(\nabla_K \ploc{K})_{K\in\mathcal M}$ does not allow to control $p$ itself);
the important characteristic of $S_K$ is that it yields a coercive bilinear form, so that  we may in fact replace \refe{deffluxhfv1} by the more general  equation
\begin{eqnarray}
\forall K\in\mathcal M\,,\;
\forall \qloc{K}\in \Xmbloc{K}\,:&&\nonumber\\
\sum_{\sigma\in\mathcal E_K} |\sigma| F_{K,\sigma}(q_K-q_\sigma)
&=&|K|\Lambda_K \nabla_K \ploc{K}\cdot\nabla_K \qloc{K}+\sum_{\sigma,\sigma'\in\mathcal E_K}
\mathbb{B}^H_{K,\sigma,\sigma'}S_{K,\sigma}(\ploc{K})S_{K,\sigma'}(\qloc{K})\nonumber\\
&=&|K|\Lambda_K \nabla_K \ploc{K}\cdot\nabla_K \qloc{K}+S_{K}(\qloc{K})^T\mathbb{B}^H_{K}S_{K}(\ploc{K}),
\label{deffluxhfv2}\end{eqnarray}
where $\mathbb{B}^H_{K,\sigma,\sigma'}$ are the entries of a  symmetric definite positive matrix $\mathbb{B}^H_K$ (the superscript $H$ refers to the Hybrid method).

\begin{defi}[Generalized Hybrid method] A Generalized Hybrid scheme for \refe{base} reads:
\begin{center}
Find $(\widetilde{p},F)\in \Xmbzero\times\mathcal F$ which satisfies \listeeq{\refe{defgradhfv},\refe{defSk},\refe{defhfv_bis},\refe{deffluxhfv2}}.
\end{center}
Its parameters are the  family of points $(x_K)_{K\in\mathcal M}$ (which are freely chosen inside each grid cell)
and the family  $(\mathbb{B}^H_K)_{K\in\mathcal M}$ of symmetric definite positive matrices.
\label{def-genhybrid}\end{defi}

\br Another presentation of the Generalized Hybrid method is possible from \refe{defhfv} and \refe{deffluxhfv2} by eliminating the fluxes: 
\begin{center}
Find $\widetilde{p}\in \Xmbzero$ such that
\end{center}
\begin{equation}
\forall \widetilde{q}\in \Xmbzero\,:\;
\sum_{K\in\mathcal M}|K|\Lambda_K \nabla_K \ploc{K}\cdot\nabla_K \qloc{K}
+\sum_{K\in\mathcal M} S_K(\qloc{K})^T
\mathbb{B}^H_K S_{K}(\ploc{K})
=\sum_{K\in\mathcal M}q_K\int_K f.
\label{defhfv2}\end{equation}
\label{rem-defhfv2}\er

\subsection{The Mixed method}\label{sec-mfv}

As in the Hybrid method,  we use the unknowns $\widetilde{p}$ in $\Xmb$ for the Mixed method (that is to say approximate values of the solution to the  equation inside the control volumes and on the edges), and fluxes unknowns in $\mathcal F$.
However,  contrary to the Hybrid method, the discrete gradient is  no longer defined from $p$, but rather from the fluxes,
using the dual version of \eqref{formmag2}, that is:
\begin{equation}
|K| \mathbf{e}=\sum_{\sigma\in\mathcal E_K} |\sigma| \mathbf{e}\cdot\n_{K,\sigma}
(\xs-x_K).
\label{formmag}
\end{equation}
For $F\in\mathcal F$, we denote $F_K=(F_{K,\sigma})_{\sigma\in\mathcal E_K}$ its restriction to the edges of the control volume $K$ and $\mathcal F_K$ is the set of such restrictions; if $F_K\in\mathcal F_K$, then $\v_K(F_K)$ is defined by
\begin{equation}\label{defvmfv}
|K|\Lambda_K \v_K(F_K)=-\sum_{\sigma\in\mathcal E_K}|\sigma|F_{K,\sigma}(\xs-x_K)
\end{equation}
where, again, $\Lambda_K$ is the mean value of $\Lambda$ on $K$, $\xs$ is the center of
gravity of $\sigma$ and $x_K$ a point chosen inside $K$.
Recalling that $F_{K,\sigma}$ is an approximation of $\frac{1}{|\sigma|}\int_\sigma -\Lambda\nabla p\cdot\n_{K,\sigma}$, Formula
\refe{formmag}  then shows that $\v_K(F_K)$ is expected to play the role of an approximation of $\nabla p$ in $K$.

The Mixed method then consists in seeking $(\widetilde{p},F)\in\Xmbzero\times \mathcal F$ (recall that $F \in \mathcal F $ satisfies \refe{conservativite}, and we impose $p_\sigma=0$ if $\sigma\subset \partial\Omega$) which satisfies the following natural discrete relation between $p$ and this discrete gradient, with a stabilization term involving the fluxes and  a positive parameter $\nu>0$:
\begin{equation}
\forall K\in\mathcal M\,,\;\forall\sigma\in\mathcal E_K\,:\;
p_\sigma-p_K=\v_K(F_K)\cdot (\xs-x_K)-\nu{\rm diam}(K)F_{K,\sigma}
\label{mfv-penforte}\end{equation}
and the flux balance:
\begin{equation}
\forall K\in\mathcal M\,:\;
\sum_{\sigma\in\mathcal E_K}|\sigma|F_{K,\sigma}=\int_K f.
\label{bilanmfv}\end{equation}

\medskip

Multiplying, for any $G_K\in\mathcal F_K$,
each equation of \refe{mfv-penforte} by $|\sigma|G_{K,\sigma}$ and summing
on $\sigma\in\mathcal E_K$, we obtain
\begin{eqnarray*}
&&\forall K\in\mathcal M\,,\;\forall G_K\in\mathcal F_K\,:\\
&&\sum_{\sigma\in\mathcal E_K}(p_K-p_\sigma)|\sigma|G_{K,\sigma}
=|K|\v_K(F_K)\cdot\Lambda_K\v_K(G_K)+\sum_{\sigma\in\mathcal E_K}
\nu {\rm diam}(K)|\sigma| F_{K,\sigma}G_{K,\sigma}.
\end{eqnarray*}
The term $\sum_{\sigma\in\mathcal E_K}
\nu {\rm diam}(K)|\sigma|F_{K,\sigma}G_{K,\sigma}$ in this equation can be considered as a strong stabilization term, in the sense that it vanishes (for all $G_K$) only if $F_{K}$ vanishes.
We modify here the Mixed method by replacing this strong stabilization by
a weaker stabilization which, as in the Hybrid method, is expected to vanish  on ``appropriate fluxes''.
To achieve this, we use the quantity
\begin{equation}
T_{K,\sigma}(F_K)=F_{K,\sigma}+\Lambda_K \v_K(F_K)\cdot\n_{K,\sigma},
\label{defpenmfv}
\end{equation}
which vanishes if $(F_{K,\sigma})_{\sigma\in\mathcal E_K}$ are the genuine fluxes of a given vector $\mathbf{e}$. 
Then, taking $\mathbb{B}^M_K$ to be a symmetric positive definite matrix, we endow $\mathcal F_K$ with the inner product
\begin{eqnarray}
\langle F_K,G_K\rangle_K&=&|K|\v_K(F_K)\cdot\Lambda_K\v_K(G_K)+\sum_{\sigma,\sigma'\in\mathcal E_K}
\mathbb{B}^M_{K,\sigma,\sigma'}T_{K,\sigma}(F_K)T_{K,\sigma'}(G_K)\nonumber\\
&=&|K|\v_K(F_K)\cdot\Lambda_K\v_K(G_K)+
T_K(G_K)^T\mathbb{B}^M_K T_{K}(F_K)
\label{defpslocalmfv}\end{eqnarray}
and the stabilized formula \refe{mfv-penforte} linking $p$ and $F$ is replaced by
\begin{equation}
\forall K\in\mathcal M\,,\;\forall G_K\in\mathcal F_K\,:\;
\langle F_K,G_K\rangle_K=\sum_{\sigma\in\mathcal E_K}(p_K-p_\sigma)|\sigma|G_{K,\sigma}.
\label{lienpFmfv}\end{equation}
We can get back a formulation more along the lines of \refe{mfv-penforte}
if we fix $\sigma\in \mathcal E_K$ and take $G_K(\sigma)\in\mathcal F_K$
defined by $G_K(\sigma)_\sigma=1$ and $G_K(\sigma)_{\sigma'}=0$ if $\sigma'\not=\sigma$:
\refe{lienpFmfv} then gives
\[
p_\sigma-p_K=-\frac{1}{|\sigma|}|K|\v_K(F_K)\cdot\Lambda_K\v_K(G_K(\sigma))
-\frac{1}{|\sigma|}T_K(G_K(\sigma))^T\mathbb{B}^M_KT_K(F_K).
\]
But $|K|\Lambda_K\v_K(G_K(\sigma))=-|\sigma|(\xs-x_K)$ and thus
\[
p_\sigma-p_K=\v_K(F_K)\cdot (\xs-x_K)-\frac{1}{|\sigma|}T_K(G_K(\sigma))^T\mathbb{B}^M_KT_K(F_K).
\]
This equation is precisely the natural discrete relation \refe{mfv-penforte} between $p$ and its gradient, in which the strong stabilization involving $F_{K,\sigma}$ has been replaced by a ``weak'' stabilization using $T_K(F_K)$.

\begin{defi}[Modified Mixed method] A Modified Mixed scheme for \eqref{base} reads:
\begin{center}
Find $(\widetilde{p},F)\in \Xmbzero\times\mathcal F$ which satisfies 
\listeeq{\refe{defvmfv},\refe{bilanmfv},\refe{defpenmfv},\refe{defpslocalmfv},\refe{lienpFmfv}}.
\end{center}
Its parameters are the  family of points $(x_K)_{K\in\mathcal M}$ (which are freely chosen inside each grid cell)
and the family  $(\mathbb{B}^M_K)_{K\in\mathcal M}$ of symmetric definite positive matrices.
\label{def-modmixed}\end{defi}

\section{Algebraic correspondence between the three methods}\label{algebraic}

We now focus on the main result of this paper, which is the following.

\begin{theorem}[Equivalence of the methods] The Generalized Mimetic, Generalized Hybrid and Mo\-dified Mixed methods (see Definitions \ref{def-genmimetic}, \ref{def-genhybrid} and \ref{def-modmixed}) are identical, in the sense that for any choice of parameters for one of these methods, there exists a choice of parameters for the other two methods which leads to the same scheme.
\label{main}\end{theorem}
The proof of this result is given in the remainder of this section, and decomposed as follows:
for comparison purposes, the  Generalized Mimetic method is first written under a hybridized form  in Section \ref{sec-hybrider_mimetic};
then, the correspondence between the Generalized Mimetic and the Modified Mixed methods is studied  in Section \ref{mimix}; 
finally, the correspondence between the Generalized Mimetic and the Generalized Hybrid methods is carried out in Section \ref{sec-mfdmfv}.

\subsection{Hybridization of the Generalized Mimetic method}\label{sec-hybrider_mimetic}

Although the edge unknowns introduced to define the Generalized Hybrid and Modified Mixed methods may be eliminated, they are in fact essential to these methods; indeed,   both methods can be hybridized and reduced to a system with unknowns $(p_\sigma)_{\sigma\in\mathcal E}$ only. 
In order to compare the methods, we therefore also introduce such edge unknowns in the Generalized Mimetic method; this is the aim of this section.

Let $E$ be a grid cell and $e\in\partial E$ be an edge.
If $e$ is an interior edge, we denote by $\widetilde{E}$  the cell on the other side of $e$ and define $G(e,E)\in X^h$ by:
$G(e,E)_E^e = 1$, $G(e,E)_{\widetilde{E}}^e=-1$ and $G(e,E)_{E'}^{e'}=0$ if $e'\not= e$.
We notice that $G(e,E) = -G(e,\widetilde{E})$ and, using $G(e,E)$ in \eqref{defmfd1}, the definitions of $\divmfd$ and of the inner products on $X^h$ and $Q^h$ give $|e|(p_E p_{\widetilde{E}}) = [F,G(e,E)]_E+[F,G(e,E)]_{\widetilde{E}} =[F,G(e,E)]_E-[F,G(e,{\widetilde{E}})]_{\widetilde{E}}$.
It is therefore natural to define $p_e$ (only depending on $e$ and not  on the grid cell $E$ such that $e\subset\partial E$) by
\begin{equation}
\forall E\in \Omega_h\,,\;\forall e\in\partial E\,:\;|e|(p_E-p_e)=[F,G(e,E)]_E.
\label{defpemfd}\end{equation}
This definition can also be applied for boundary edges $e\subset \partial\Omega$,
in which case it gives $p_e=0$ (thanks to \eqref{defmfd1}).

We thus extend $p\in Q^h$ into an element $\widetilde{p}\in \Xmbzero$ having values inside the cells and on the edges of the mesh. 
Denoting, as in Section \ref{sec-mfv}, $\mathcal F_E$ the space of restrictions $G_E$ to the  edges of $E$ of elements $G\in X^h$ and writing any $G_E\in\mathcal F_E$ as a linear 
combination of $(G(E,e))_{e\in\partial E}$, it is easy to see from \eqref{defpemfd} that the Generalized Mimetic method \listeeq{\refe{defmfd1},\refe{defmfd2}} is equivalent to: 
find  $(\widetilde{p},F)\in\Xmbzero\times X^h$ such that
\begin{equation}
\forall E\in\Omega_h\,,\;\forall G_E\in\mathcal F_E\,:\;
[F_E,G_E]_E=\sum_{e\in\partial E}|e|(p_E-p_e)G_E^e
\label{mfdhyb1}\end{equation}
and
\begin{equation}
\forall E\in\Omega_h\,:\;\sum_{e\in\partial E}|e|F_E^e=\int_E f.
\label{mfdhyb2}
\end{equation}

\subsection{Proof of the correspondence Generalized Mimetic $\leftrightarrow$ Modified Mixed}\label{mimix}

The simplest comparison is probably to be found between the Modified Mixed  and Generalized Mimetic methods.
Indeed, we see from \listeeq{\refe{bilanmfv},\refe{lienpFmfv}} and \listeeq{\refe{mfdhyb1},\refe{mfdhyb2}}
that both methods are identical provided that, for any grid cell
$K=E$, the local inner products defined by \eqref{defpslocalmfv} and \eqref{choixmatmfd} are equal: $\langle\cdot,\cdot\rangle_K=[\cdot,\cdot]_E$,. 
We therefore have to study these local inner products and understand whether they can be identical (recall that there is some latitude in the choice of both inner products).

Let us start with the term  $|K| \v_K(F_K)\cdot\Lambda_K \v_K(G_K)$ in the definition of  $\langle\cdot,\cdot\rangle_K$.
Thanks to \refe{defvmfv},
\[
|K|\v_K(F_K)\cdot\Lambda_K\v_K(G_K)=\frac{1}{|K|}\left(\sum_{\sigma\in \mathcal E_K}
\Lambda_K^{-1}\big(|\sigma|(\xs-x_K)\big)F_{K,\sigma}\right)
\cdot\left(\sum_{\sigma\in \mathcal E_K} |\sigma|(\xs-x_K)G_{K,\sigma}\right).
\]
But, with the definition \refe{defnewR} of $\mathbb{R}_E$,
$\sum_{\sigma\in \mathcal E_K}|\sigma|(\xs-x_K)F_{K,\sigma}=\mathbb{R}_E^TF_K$
and thus
\[
|K|\v_K(F_K)\cdot\Lambda_K\v_K(G_K)=\frac{1}{|K|}\left(\Lambda_K^{-1}\mathbb{R}_E^TF_K\right)
\cdot\left(\mathbb{R}_E^T G_{K}\right)
=G_{K}^T \left(\frac{1}{|K|}\mathbb{R}_E\Lambda_K^{-1}\mathbb{R}_E^T\right) F_K.
\]
Hence, the term  $|K|\v_K(\cdot)\cdot\Lambda_K\v_K(\cdot)$ in the definition of  $\langle \cdot,\cdot\rangle_K$ is identical to the first term
$\frac{1}{|E|}\mathbb{R}_E\Lambda_E^{-1}\mathbb{R}_E^T$ in the definition of the matrix $\mathbb{M}_E$ of $[\cdot,\cdot]_E$ (see \refe{choixmatmfd}). 
Therefore, in order to prove that $\langle\cdot,\cdot\rangle_K=[\cdot,\cdot]_E$, it only remains to prove that, for appropriate choices of $\mathbb{C}_E$, $\mathbb{U}_E$ and $\mathbb{B}_K^M$, we have for any $(F_K,G_K) \in \mathcal F_K^2$:
\begin{equation}
T_K(G_K)^T\mathbb{B}^M_KT_K(F_K)=G_K^T\mathbb{C}_E\mathbb{U}_E\mathbb{C}_E^T F_K
\label{eg-mfdmfv}\end{equation}
(see \refe{defpslocalmfv} and \refe{choixmatmfd}); in fact, this is  the consequence
of Lemma \ref{lem-clef} in the appendix and of the following lemma.

\bl The mappings $T_K:\R^{k_E}\to \R^{k_E}$ and $\mathbb{C}_E^T:\R^{k_E}\to \R^{k_E-d}$
have the same kernel.
\label{lem-noy}\el

\bpr{ of Lemma \ref{lem-noy}}

We first prove that:
\begin{itemize}
\item[i)] ${\rm Im}(\mathbb{N}_E)\subset {\rm Ker}(T_K)$
i.e. $T_{K,\sigma}((\mathbb{N}_E)_j)=0$ for all $\sigma\in\mathcal E_K$ and all $j=1,\ldots,d$
(which also amounts to the fact that the lines of $T_K$ are orthogonal to the vectors $(\mathbb{N}_E)_j$).
\item[ii)] ${\rm dim}({\rm Im}(T_K))=k_E-d$, and therefore ${\rm dim}({\rm Ker}(T_K))=d$.
\end{itemize}

Item i) follows from \refe{defvmfv} and \refe{formmag}:
we have $\Lambda_K\v_K((\mathbb{N}_E)_j)=-\frac{1}{|K|}\sum_{\sigma\in\mathcal E_K}
|\sigma| (\Lambda_K)_j\cdot\n_{K,\sigma}(\xs-x_K)=-(\Lambda_K)_j$
and thus $T_{K,\sigma}((\mathbb{N}_E)_j)=(\Lambda_K)_j\cdot\n_{K,\sigma}-
(\Lambda_K)_j\cdot\n_{K,\sigma}=0$.

In order to obtain Item ii), we first remark that $k_E-d$ is an upper bound of the rank of $T_K$ since the lines
of $T_K$ are in the orthogonal space of the independent vectors $((\mathbb{N}_E)_j)_{j=1,\ldots,d}$
(\footnote{Let us notice that the independence of $((\mathbb{N}_E)_j)_{j=1,\ldots,d}$ can be deduced from the independence
of the columns of $\Lambda_K$: thanks to \refe{formmag}, any non-trivial
linear combination of the $(\mathbb{N}_E)_j$ gives a non-trivial combination of the columns of $\Lambda_K$.}). 
Moreover, \refe{defvmfv} shows that the rank of
$\v_K:\R^{k_E}\to \R^{k_E}$ is the rank of the family $(\xs-x_K)_{\sigma\in\mathcal E_K}$,
that is to say $d$, and the kernel of $\v_K$ therefore has dimension $k_E-d$;
since $T_K={\rm Id}$ on this kernel, we conclude that the rank of $T_K$ is at least $k_E-d$,
which proves ii).

These properties show that  ${\rm Ker}(T_K)={\rm Im}(\mathbb{N}_E)=({\rm Im}(\mathbb{C}_E))^\bot=\ker(\mathbb{C}_E^T)$, and the proof is complete. \epr

\medskip

The comparison between $T_K(G_K)^T\mathbb{B}^M_KT_K(F_K)$ and $G_K^T\mathbb{C}_E\mathbb{U}_E\mathbb{C}_E^TF_K$
is now straightforward. Indeed, let $(\mathbb{C}_E,\mathbb{U}_E)$ be any pair chosen to
construct the Generalized Mimetic method; applying Lemma \ref{lem-clef}, thanks to Lemma \ref{lem-noy},
with $A=\mathbb{C}_E^T$, $B=T_K$ and $\{\cdot,\cdot\}_{\R^{k_E-d}}$ the inner product on $\R^{k_E-d}$
corresponding to $\mathbb{U}_E$, we obtain an inner product $\{\cdot,\cdot\}_{\R^{k_E}}$ on $\R^{k_E}$
such that $\{T_K(F_K),T_K(G_K)\}_{\R^{k_E}}=\{\mathbb{C}_E^T F_K,\mathbb{C}_E^T G_K\}_{\R^{k_E-d}}=
G_K^T \mathbb{C}_E\mathbb{U}_E\mathbb{C}_E^T F_K$;
this exactly means that, if we define $\mathbb{B}_K^M$ as the matrix of $\{\cdot,\cdot\}_{\R^{k_E}}$,
\refe{eg-mfdmfv} holds.
Similarly, inverting the role of $\mathbb{C}_E^T$ and $T_K$ when applying Lemma \ref{lem-clef},
for any $\mathbb{B}^M_K$ used in the Modified Mixed method we can find
$\mathbb{U}_E$ satisfying \refe{eg-mfdmfv} and the proof of the correspondence between the
Generalized Mimetic and Modified Mixed methods is concluded.

\subsection{Proof of the correspondence Generalized Mimetic $\leftrightarrow$ Generalized Hybrid}\label{sec-mfdmfv}

To compare the Generalized Mimetic and Generalized Hybrid methods, we use a result of \cite{bre-05-fam,man-08-mim} which states that the inverse of the matrix $\mathbb{M}_E$ in \refe{choixmatmfd} can be written
\begin{equation}
\mathbb{M}_E^{-1}=\mathbb{W}_E=\frac{1}{|E|}\mathbb{N}_E\Lambda_E^{-1}\mathbb{N}_E^T+\mathbb{D}_E\widetilde{\mathbb{U}}_E
\mathbb{D}_E^T
\label{matW}\end{equation}
where $\mathbb{D}_E$ is  a $k_E\times (k_E-d)$ matrix such that 
\begin{equation}
{\rm Im} (\mathbb{D}_E) = ({\rm Im} (\mathbb{R}_E) )^\perp \label{dere}
\end{equation} 
and $\widetilde{\mathbb{U}}_E$ is a symmetric definite positive $(k_E-d)\times (k_E-d)$ matrix
(note that the proof in \cite{bre-05-fam,man-08-mim} assumes $x_E$ to be the center of gravity of $E$, i.e.   $\mathbb{M}_E$ to satisfy \refe{choixmatmfd-initial},
but that it remains valid for any choice of $x_E$, i.e. for any matrix $\mathbb{M}_E$ satisfying \refe{choixmatmfd}).
This result is to be understood in the following sense: for any $(\mathbb{C}_E,\mathbb{U}_E)$ used
to construct $\mathbb{M}_E$ by \refe{choixmatmfd}, there exists $(\mathbb{D}_E,\widetilde{\mathbb{U}}_E)$
such that $\mathbb{W}_E$ defined by \refe{matW} satisfies $\mathbb{W}_E=\mathbb{M}_E^{-1}$,
and \emph{vice-versa}.

\medskip

For $\ploc{E}=(p_E,(p_e)_{e\in\partial E})\in \Xmbloc{E}$, we denote
by $\mathcal P_E$ the vector in $\R^{k_E}$ with components $(|e|(p_E-p_e))_{e\in\partial E}$.
The relation \refe{mfdhyb1} can be re-written $\mathbb{M}_EF_E=\mathcal P_E$,
that is to say $F_E=\mathbb{W}_E \mathcal P_E$, which is also equivalent,
taking the inner product in $\R^{k_E}$ with an arbitrary $\mathcal Q_E$ (built from
a $\qloc{E}\in\Xmbloc{E}$), to
\begin{equation}
\forall E\in\Omega_h\,,\;\forall \qloc{E}\in\Xmbloc{E}\,:\;
\sum_{e\in\partial E}|e|(q_E-q_e)F_E^e=\mathcal Q_E^T\mathbb{W}_{E}\mathcal P_E
\label{defmfdhyb3}\end{equation}
The Generalized Mimetic method \listeeq{\refe{mfdhyb1},\refe{mfdhyb2}} is therefore identical to the
Generalized Hybrid method \listeeq{\refe{defhfv_bis},\refe{deffluxhfv2}} provided that, for
all $E=K\in \Omega_h$,
\begin{equation}
\forall (\ploc{E},\qloc{E})\in \Xmbloc{E}\,:\;
\mathcal Q_E^T\mathbb{W}_{E}\mathcal P_E
=|K|\Lambda_K\nabla_K \ploc{K}\cdot\nabla_K \qloc{K}+S_K(\qloc{K})^T\mathbb{B}_K^H S_K(\ploc{K}).
\label{aprouver}\end{equation}
As in the comparison between the Generalized Mimetic and Modified Mixed methods, the proof of \refe{aprouver}
is obtained from the separate study of each term of \refe{matW}.

\medskip

First, by  the definition \eqref{mf2} of the matrix $\mathbb{N}_E$ and  the definition \eqref{defgradhfv} of $\nabla_K \ploc{K}$,
we have $(\mathbb{N}_E^T\mathcal P_E)_j=\sum_{e\in\partial E}|e|(\Lambda_E)_j\cdot\n^e_E (p_E-p_e)=
-|K|(\Lambda_K)_j\cdot \nabla_K \ploc{K}$ for all $j=1,\ldots,d$, that is to say,
by symmetry of $\Lambda$, $\mathbb{N}_E^T\mathcal P_E=-|K|\Lambda_K\nabla_K\ploc{K}$.
Hence,
\[
\mathcal Q_E^T\left(\frac{1}{|E|}\mathbb{N}_E\Lambda_K^{-1}\mathbb{N}_E^T\right)\mathcal P_E
=|K|(\Lambda_K\nabla_K\qloc{K})^T \Lambda_K^{-1} (\Lambda_K\nabla_K\ploc{K})
=|K|\Lambda_K\nabla_K \ploc{K}\cdot\nabla_K \qloc{K}.
\]
The first part of the right-hand side in \refe{aprouver} thus corresponds to the
first part of $\mathbb{W}_E$
in \refe{matW}, and it remains to prove that (with appropriate choices of $\mathbb{D}_E$,
$\widetilde{\mathbb{U}}_E$ and $\mathbb{B}_K^H$), for
all $(\ploc{E},\qloc{E})\in \Xmbloc{E}$,
$\mathcal Q_E^T\mathbb{D}_E\widetilde{\mathbb{U}}_E\mathbb{D}_E^T\mathcal P_E=
S_K(\qloc{K})^T\mathbb{B}_K^H S_K(\ploc{K})$.
Let us notice that, defining $L_K:\R^{k_E}\to \R^{k_E}$ by
\begin{equation}
L_K(V)=\left(\frac{1}{|\sigma|}V_\sigma-D_K V \cdot (\xs-x_K)\right)^T_{\sigma\in\mathcal E_K}
\quad\mbox{ with }\quad D_K V = \frac{1}{|K|}\sum_{\sigma\in\mathcal E_K}V_\sigma\n_{K,\sigma},
\label{deflk}\end{equation}
this boils down (letting $V=\mathcal P_E$ and $V'=\mathcal Q_E$) to proving that
\begin{equation}
\forall (V,V')\in \R^{k_E}\,:\; (V')^T\mathbb{D}_E\widetilde{\mathbb{U}}_E\mathbb{D}_E^T V
=L_K(V')^T\mathbb{B}_K^H L_K(V).
\label{aprouver2}\end{equation}
As previously, this will be a consequence of Lemma \ref{lem-clef} and of the following result.

\bl The mappings $L_K:\R^{k_E}\to \R^{k_E}$ and $\mathbb{D}_E^T:\R^{k_E}\to \R^{k_E-d}$ have the same kernel.
\label{lem-noy2}\el

\bpr{ of Lemma \ref{lem-noy2}}
From \eqref{dere}, we get that ${\rm Ker} (\mathbb{D}_E^T) = {\rm Im}  (\mathbb R_E)$. 
Hence it remains to prove that ${\rm Im}  (\mathbb R_E) = {\rm Ker}(L_K)$

Let us first prove that ${\rm Im}  (\mathbb R_E )\subset {\rm Ker}(L_K)$.
The $j$-th column of $\mathbb{R}_E$ is $(\mathbb{R}_E)_j=(|\sigma|(\xs^j-x_K^j))^T_{\sigma\in\mathcal E_K}$
(the superscript $j$ denotes the $j$-th coordinate of points in $\R^d$).
Thus, for any $\mathbf{e}\in \R^d$, by \refe{formmag}, 
$$D_K(\mathbb{R}_E)_j\cdot \mathbf{e}=\frac{1}{|K|}\sum_{\sigma\in\mathcal E_K} |\sigma|\mathbf{e}\cdot\n_{K,\sigma} (\xs^j-x_K^j)= \mathbf{e}^j,$$ which means that $D_K(\mathbb{R}_E)_j$ is the $j$-th vector of the canonical basis of $\R^d$. 
We then have $D_K(\mathbb{R}_E)_j\cdot(\xs-x_K)=\xs^j-x_K^j$ and thus 
$$(L_{K}((\mathbb{R}_E)_j))_\sigma=\xs^j-x_K^j-D_K(\mathbb{R}_E)_j\cdot(\xs-x_K)=0;$$
this proves that the columns
of $\mathbb{R}_E$ are in the kernel of $L_K$, that is ${\rm Im}  (\mathbb R_E) \subset  {\rm Ker}  (L_K)$.

Next we notice that the rank of the mapping $D_K$:  $V\in\R^{k_E}\mapsto D_K V\in\R^{d}$ (i.e. the rank of the family 
$(\n_{K,\sigma})_{\sigma\in\mathcal E_K}$) is $d$ and that the mapping  $L_K$ is one-to-one on ${\rm Ker} (D_K)$. 
Hence $\mathrm{dim}(\mathrm{Im}(L_K))\ge {\rm dim}({\rm Ker}(D_K))=k_E-d$, and therefore
$\mathrm{dim}({\rm Ker}(L_K)) \le  d$.
But since $\mathrm{Im}(\mathbb R_E) \subset \mathrm{Ker}(L_K)$ and ${\rm dim}({\rm Im}(\mathbb{R}_E)=d$
(the rank of the rows of $\mathbb{R}_E$),
we thus get that $\mathrm{ Im}(\mathbb R_E) = \mathrm {Ker}(L_K)$. \epr

\medskip

{}From Lemmas \ref{lem-clef} and \ref{lem-noy2}, we deduce as in Section
\ref{mimix} that for any choice of $(\mathbb{D}_E,\widetilde{\mathbb{U}}_E)$
there corresponds a choice of $\mathbb{B}_K^H$ such that \refe{aprouver2}
holds, and \emph{vice-versa}, which concludes the proof that the Generalized Mimetic
method is identical to the Generalized Hybrid method.

\br These proofs and Remark \ref{rem-explicite0} show
that one can explicitly compute the parameters of one method which
gives back the parameters of another method. Notice also that the algebraic computations
required in Remark \ref{rem-explicite0} to obtain these parameters are
made in spaces with small dimensions; the cost of their practical
computations is therefore very low. However,
to implement the methods, it is not necessary to understand which $(\mathbb{C}_E,\mathbb{U}_E)$
or $(\mathbb{D}_E,\widetilde{\mathbb{U}}_E)$ corresponds to which $\mathbb{B}^M_K$
or $\mathbb{B}^H_K$, since the only useful quantities are $\mathbb{C}_E\mathbb{U}_E\mathbb{C}_E^T$,
$\mathbb{D}_E\widetilde{\mathbb{U}}_E\mathbb{D}_E^T$,
$T_K(\cdot)^T\mathbb{B}^M_KT_K(\cdot)$ and $L_K(\cdot)^T\mathbb{B}^H_K
L_K(\cdot)$, and the relations between these quantities are
trivial (see \refe{eg-mfdmfv} and \refe{aprouver2}).
\label{rem-explicite}\er

\section{Convergence and error estimates} \label{sec-convergence}

We showed in Section \ref{algebraic} that the three families of schemes which we called Generalized Mimetic, Generalized Hybrid FV and Modified Mixed FV are in fact one family of schemes, which we call the HMMF (Hybrid Mimetic Mixed Family) for short in the remainder of the paper.
We now give  convergence  and error estimate results for the HMMF method. 

\subsection{Convergence with no regularity}

In this section, we are interested in convergence results which hold without any other regularity assumption on the data than those stated in Section \ref{sec-intro}. 
We therefore consider that $\Lambda$ is only bounded and uniformly elliptic (not necessarily Lipschitz continuous or even piecewise continuous), that $f\in L^2(\Omega)$ and that the solution to \eqref{base} only belongs to $H^1_0(\Omega)$ (and not necessarily to $H^2(\Omega)$).

We study how existing results, previously established  for the Hybrid scheme, can be extended to the HMMF.
In \cite{eym-08-dis}, we proved the $L^2$ convergence of the ``standard'' Hybrid method for a family of partitions of $\Omega$ such that any cell $K$ is star-shaped with respect to a point $x_K$ and such that there exists $\theta>0$ satisfying, for any partition of the family,
\begin{equation}
\max\left(\max_{\substack{\sigma\in\mathcal{E}_{\mathrm{int}} \\ K,L \in  \mathcal{M}_\sigma} }\frac {d_{K,\sigma}} {d_{L,\sigma}}, \max_{\substack{K \in \mathcal{M} \\ \sigma \in\mathcal{E}_K} }\dfrac{{\rm diam}(K)}{d_{K,\sigma}} \right) \le \theta,
\label{reghyb}
\end{equation}
where $\mathcal{E}_{\mathrm{int}}$ denotes the set of internal edges of the mesh and $\mathcal{M}_\sigma$ the set of cells to which $\sigma$ is an edge.
We now show how the convergence of the HMMF  may be deduced from an easy extension of   \cite[Theorem 4.1, Lemma 4.4]{eym-08-dis} provided that:
\begin{itemize}
\item in the Generalized Mimetic formulation 
\listeeq{\refe{defdivmfd},\refe{localscal},\refe{defmfd1},\refe{defmfd2},\refe{choixmatmfd},\refe{defnewR}},
\begin{equation}\left.\begin{array}{c}
\hbox{there exist  }s_*>0\hbox{ and }S_*>0,\hbox{ independent of the mesh, such that},\\
\hbox{for all cell }K\hbox{ and all }V=(V_\sigma)_{\sigma\in\mathcal E_K},\\
\dsp s_* \sum_{\sigma\in\mathcal E_K} |\sigma| d_{K,\sigma} (V_\sigma)^2 \le
V^T\mathbb{M}_K V
\le S_* \sum_{\sigma\in\mathcal E_K} |\sigma| d_{K,\sigma}  (V_\sigma)^2,
\end{array}\right\}\label{S1hyb1}\end{equation}
\item in the Generalized Hybrid formulation
\listeeq{\refe{defgradhfv},\refe{defSk},\refe{defhfv_bis},\refe{deffluxhfv2}},
using the notation \refe{deflk},
\begin{equation}\left.\begin{array}{c}
\hbox{there exist  }\bar s_*>0\hbox{ and }\bar S_*>0,\hbox{ independent of the mesh, such that},\\
\hbox{for all cell }K\hbox{ and all }V=(V_\sigma)_{\sigma\in\mathcal E_K},\\
\bar s_* \sum_{\sigma \in \mathcal{E}_K} \dfrac{\vert\sigma\vert}{d_{K,\sigma}} (L_K(V))_\sigma^2\le 
L_K(V)^T\mathbb{B}_K^H L_K(V) \le \bar S_* \sum_{\sigma \in \mathcal{E}_K} \dfrac{\vert\sigma\vert}{d_{K,\sigma}} 
(L_K(V))_\sigma^2, \label{c1c2}
\end{array}\right\}\end{equation}
\item in the Modified Mixed formulation
\listeeq{\refe{defvmfv},\refe{bilanmfv},\refe{defpenmfv},\refe{defpslocalmfv},\refe{lienpFmfv}},
\begin{equation}\left.\begin{array}{c}
\hbox{there exist  }\tilde s_*>0\hbox{ and }\tilde S_*>0,\hbox{ independent of the mesh, such that},\\
\hbox{for all cell }K\hbox{ and all }V=(V_\sigma)_{\sigma\in\mathcal E_K},\\
\tilde s_*\sum_{\sigma \in \mathcal{E}_K} |\sigma|d_{K,\sigma} (T_{K,\sigma}(V))^2\le 
T_K(V)^T\mathbb{B}_K^M T_K(V)\le \tilde S_*\sum_{\sigma \in \mathcal{E}_K} |\sigma| d_{K,\sigma}(T_{K,\sigma}(V))^2.
\end{array}\right\}\label{cond-matmfv}\end{equation}
\end{itemize}
The three conditions \eqref{S1hyb1}, \eqref{c1c2} and \eqref{cond-matmfv} are in fact equivalent (this is stated in the next theorem), and one only has to check the condition corresponding to the chosen framework.

\begin{theorem}[Convergence of the HMMF method]\label{convhfvgen}
Assume that $\Lambda:\Omega\to M_d(\R)$  is bounded measurable symmetric and uniformly elliptic, that $f\in L^2(\Omega)$ and that the solution to \eqref{base} belongs to $H^1_0(\Omega)$.
Let $\theta >0$ be given.
Consider a family of polygonal meshes of $\Omega$ such that any cell $K$ is star-shaped with respect to a point $x_K$, and satisfying \eqref{reghyb}.
Then the three conditions \eqref{S1hyb1}, \eqref{c1c2} and \eqref{cond-matmfv} are equivalent. 
Moreover, if for any mesh of the family we choose a HMMF scheme such that one of the conditions \eqref{S1hyb1}, \eqref{c1c2} or \eqref{cond-matmfv} holds,
then, the family of corresponding approximate solutions converges in $L^2(\Omega)$ to the unique solution of \eqref{base} as the mesh size tends to 0.
\end{theorem}

\bpr{ of Theorem \ref{convhfvgen}.}

Let us first prove the equivalence between \eqref{S1hyb1} and \eqref{c1c2}, assuming \eqref{S1hyb1} to begin with.
Using \eqref{aprouver2}, we get that, for all $V$,
\[
V^T\mathbb{D}_K\widetilde{\mathbb{U}}_K\mathbb{D}_K^T V
=L_K(V)^T\mathbb{B}_K^H L_K(V).
\]
Let us apply the above relation to  $\widetilde L_K(V)$ defined by $\widetilde L_K(V)_\sigma = |\sigma|L_K(V)_\sigma$. From \eqref{formmag}
(see also \eqref{formmag2}) and recalling the operator $D_K$ defined in \eqref{deflk}, we easily get that $D_K \widetilde L_K(V) = 0$, which provides $L_K(\widetilde L_K(V))= L_K(V)$. Hence we get
\[
(\widetilde L_K(V))^T\mathbb{D}_K\widetilde{\mathbb{U}}_K\mathbb{D}_K^T \widetilde L_K(V)
=L_K(V)^T\mathbb{B}_K^H L_K(V).
\]
We then remark that $\mathbb{N}_K^T \widetilde L_K(V) = 0$ (once again from \eqref{formmag})
and therefore, using \eqref{matW},
\begin{equation}\label{coerB1}
L_K(V)^T\mathbb{B}_K^H L_K(V) = (\widetilde L_K(V))^T \mathbb{W}_K\widetilde L_K(V).
\end{equation}

Let $\mathbb{I}_K$ be the diagonal matrix
$\mathbb{I}_K={\rm diag}((\sqrt{|\sigma|d_{K,\sigma}})_{\sigma\in\mathcal E_K}))$.
Condition \eqref{S1hyb1} applied to $G_K=\mathbb{I}_KF_K$ gives, for all $G_K\in\mathcal F_K$,
\[
s_* G_K^TG_K\le G_K^T\mathbb{I}_K^{-1}\mathbb{M}_K\mathbb{I}_K^{-1}G_K\le S_*G_K^TG_K.
\]
This shows that the eigenvalues of $\mathbb{I}_K^{-1}\mathbb{M}_K\mathbb{I}_K^{-1}$ are
in $[s_*,S_*]$, and thus that the eigenvalues of 
$\mathbb{I}_K\mathbb{M}_K^{-1}\mathbb{I}_K=
\mathbb{I}_K\mathbb{W}_K\mathbb{I}_K$
belong to $[\frac{1}{S_*},\frac{1}{s_*}]$. Translating this
into bounds on $(\mathbb{I}_K^{-1}\widetilde{L}_K(V))^T\mathbb{I}_K\mathbb{W}_K\mathbb{I}_K
(\mathbb{I}_K^{-1}\widetilde{L}_K(V))$, we deduce
\[
\frac{1}{S_*}\sum_{\sigma\in\mathcal E_K} \frac 1 {|\sigma| d_{K,\sigma}} (\widetilde L_K(V)_\sigma)^2\le (\widetilde L_K(V))^T \mathbb{W}_K\widetilde L_K(V)\le
\frac{1}{s_*}\sum_{\sigma\in\mathcal E_K} \frac 1 {|\sigma| d_{K,\sigma}} (\widetilde L_K(V)_\sigma)^2,
\]
and \eqref{c1c2} follows from \eqref{coerB1} (with $\bar s_*=\frac{1}{S_*}$ and
$\bar S_*=\frac{1}{s_*}$). Reciprocally, from \eqref{c1c2},
following the proof of \cite[Lemma 4.4]{eym-08-dis} and setting $V_\sigma = |\sigma|(p_K - p_\sigma)$ in that proof, we get the existence of $c_1>0$ and $c_2>0$, independent of the mesh, such that
\begin{equation}
c_1 \sum_{\sigma \in \mathcal{E}_K} \frac 1 {|\sigma| d_{K,\sigma}} V_\sigma^2\le 
V^T\mathbb{W}_K V \le c_2 \sum_{\sigma \in \mathcal{E}_K} \frac 1 {|\sigma| d_{K,\sigma}}
V_\sigma^2.
\label{reciproc}\end{equation}
Using $\mathbb{I}_K$ as before, we then get \eqref{S1hyb1} with $s_* = \frac{1}{c_2}$ and 
$S_* = \frac{1}{c_1}$.

\medskip

Let us turn to the proof of the equivalence between \eqref{S1hyb1} and \eqref{cond-matmfv},
beginning by assuming that \eqref{S1hyb1} holds.
Using \eqref{formmag}, one has 
$T_K((\Lambda_K\mathbf{v}_K(F_K)\cdot\n_{K,\sigma})_{\sigma\in\mathcal E_K})=0$,
and thus $T_K(T_K(F_K))=T_K(F_K)$; hence, noting that $\mathbb{R}_K^TT_K(F_K)=0$
(once again thanks to \eqref{formmag}) and remembering \eqref{choixmatmfd}, \eqref{eg-mfdmfv},
we see that \eqref{S1hyb1} applied to $V=T_K(F_K)$ directly gives \eqref{cond-matmfv}. The 
reciprocal property follows, in a similar way as \eqref{reciproc},
from a simple adaptation of classical
Mixed FV manipulations (used for example in the proof of \emph{a priori} estimates
on the approximate solution), see \cite{dro-06-mix}.

\medskip

Note also that obtaining \eqref{S1hyb1} from \eqref{c1c2} or
\eqref{cond-matmfv} is very similar to \cite[Theorem 3.6]{bre-05-fam}.

\medskip

We can now conclude the proof of the theorem, that is to say establish the
convergence of the approximate solutions: using \eqref{c1c2}, this convergence  is a direct
consequence of \cite[Theorem 4.1]{eym-08-dis} with a straightforward adaptation of the proof of
\cite[Lemma 4.4]{eym-08-dis}. The convergence could also be obtained,
using \eqref{cond-matmfv}, by an easy adaptation of the techniques of
proof used in the standard Mixed setting
(see \cite{dro-06-mix}). \epr

\medskip

Mimetic schemes are usually studied under a condition on the local inner products usually called ($\mathbf{S1}$)
and which reads \cite{bre-05-con}:
there exist $s_*>0$, $S_*>0$ independent of the mesh in the chosen family
such that
\begin{equation}
\forall E\in\Omega_h\,,\;\forall G\in X^h\,:\;
s_* \sum_{i=1}^{k_E}|E|(G_E^{e_i})^2\le [G,G]_E\le
S_* \sum_{i=1}^{k_E}|E|(G_E^{e_i})^2.
\label{S1}\end{equation}
The mesh regularity assumptions in \cite{bre-05-con,bre-05-fam} also entail
the existence of $\ctel{C3}$, independent of the mesh, such that
\begin{equation}\label{regulequiv}
\begin{array}{r}
\displaystyle\cter{C3} {\rm diam}(K)^{d-1} \le |\sigma|, \ \forall\sigma\in\mathcal E_K,\\
\displaystyle\cter{C3} {\rm diam}(K) \le d_{K,\sigma}, \ \forall\sigma\in\mathcal E_K.
\end{array}
\end{equation}
Under such mesh assumptions and since $|\sigma|\le 2^{d-1}{\rm diam}(L)$ whenever 
$\sigma\in\mathcal E_L$, it is easy to see that there exists $\theta$ only depending on 
$\cter{C3}$ such that \eqref{reghyb} holds; still using \eqref{regulequiv}, we see that
the quantities $|\sigma|d_{K,\sigma}$ are of the same order as $|K|$ and thus
that \eqref{S1} implies \eqref{S1hyb1} (with possibly different $s_*$ and $S_*$).
We can therefore apply Theorem \ref{convhfvgen} to deduce the following convergence
result under the usual assumptions in the mimetic literature
(except for the regularity assumptions on the data).

\begin{cor}[Convergence under the usual mimetic assumptions]
\label{corhfvmim}
We assume that $\Lambda:\Omega\to M_d(\R)$ 
is bounded measurable symmetric and uniformly elliptic, that
$f\in L^2(\Omega)$ and that the solution to
\eqref{base} is in $H^1_0(\Omega)$.
Consider a family of polygonal meshes of $\Omega$ such that any cell $E$ is star-shaped with 
respect to a point $x_E$, and such that
\eqref{regulequiv} holds (with $\cter{C3}$ independent of the mesh
in the family). We choose local inner products satisfying
Condition {\rm ({\bf S1})} of \cite{bre-05-con}
(i.e. \eqref{S1} with $s_*$, $S_*$ not depending on the mesh)
and \eqref{GS2}. Then, the family of approximate solutions given by the corresponding
Generalized Mimetic schemes converges in $L^2(\Omega)$ to the solution of
\eqref{base} as the mesh size tends to 0.

In particular, taking  $x_E$ as the center of gravity $\bar x_E$ of $E$,
\eqref{GS2} reduces to Condition {\rm ($\mathbf{S2}$)} of \cite{bre-05-con}
(that is to say \eqref{S2}), and the family of approximate solutions given by the
corresponding ``standard'' Mimetic schemes converges in $L^2(\Omega)$ to the
solution of \eqref{base} as the mesh size tends to 0.
\end{cor}

\begin{remark}[Compactness and convergence of a gradient]
The proofs of the above convergence results rely on the use of the
Kolmogorov compactness theorem on the family of approximate solutions, 
see \cite[Lemma 3.3]{dro-06-mix} or \cite[Section 5.2]{eym-08-dis}.
In fact, this latter study shows that  for $p \in [1,+\infty),$ any family of piecewise constant 
functions that is bounded in an adequate discrete mesh-dependent $W^{1,p}$ norm converges
in $L^p(\Omega)$ (up to a subsequence) to a function of $W^{1,p}_0(\Omega)$.
The regularity of the limit is shown thanks to the weak convergence of a discrete gradient to the 
gradient of the limit \cite[Proof of Lemma 5.7, Step 1]{eym-08-dis}. 
Note that this compactness result and the construction of this weakly converging gradient are 
completely independent of the numerical scheme, which is only used to obtain the 
discrete $W^{1,p}$  estimate and the fact that the limit is indeed the solution of \eqref{base}.

Moreover, it is possible, from the approximate solutions given by the HMMF schemes,
to reconstruct a gradient (in a similar way as \cite[(22)-(26)]{eym-08-dis})
which strongly converges in $L^2(\Omega)$ to the gradient of the solution of \eqref{base},
see \cite[Theorem 4.1]{eym-08-dis}. One can also directly prove
that the gradient $\mathbf{v}_K$, defined from the fluxes by \eqref{defvmfv},
is already a strongly convergent gradient in $L^2(\Omega)$,
see \cite{dro-06-mix}; in fact, the strong convergence of this
gradient is valid in a more general
framework than the one of the methods presented here, see \cite{agel-08-Gscheme}.
\end{remark}

\subsection{Order 1 error estimate}

Let us first consider the original Mimetic Finite Difference method \eqref{defdivmfd}--\eqref{mf3}.
This method is consistent in the sense that it satisfies the condition ({\bf S2}) of \cite{bre-05-con}.
Under the assumptions that:
\begin{itemize}
\item Condition ({\bf M1}) of \cite{bre-05-con} on the domain $\Omega$ (namely $\Omega$ is polyhedral and its  
boundary is Lipschitz continuous) and Conditions ({\bf M2})-({\bf M6}) of \cite{bre-05-con} (corresponding to 
({\bf M1})-({\bf M4}) in \cite{bre-05-fam}) hold;
\item The stability condition ({\bf S1}) of \cite{bre-05-con,bre-05-fam} holds; this condition concerns the 
eigenvalues of the matrix $\mathbb{U}_E$ in \eqref{choixmatmfd-initial}, and are connected with the discretization,
\item $\Lambda \in W^{1,\infty}(\Omega)^{d\times d}$ and $p \in H^2(\Omega)$,
\end{itemize}
then Theorem 5.2 in \cite{bre-05-con} gives an order 1 error estimate on the fluxes, for an adequate discrete norm; moreover, if $\Omega$ is convex and the right-hand side belongs to $H^1(\Omega)$, an order 1 error
estimate on $p$ in the $L^2$ norm is also established.

By the equivalence theorem (Theorem \ref{main}), this result yields similar error estimates for the Generalized Hybrid
and Modified Mixed methods in the case where the point $x_K$ is taken as the center of gravity $\bar x_K$ of $K$.

\medskip

In the case of the original Hybrid method \eqref{defhfv}, an order 1 error estimate is proved in \cite{eym-08-dis}  
only under the assumption of an homogeneous isotropic medium, that is $\Lambda = \rm{Id}$, and in the case where the 
solution to \eqref{base} also belongs to $C^2(\overline{\Omega})$ (but with no convexity assumption on the domain 
$\Omega$). 
As stated in Remark 4.2 of  \cite{eym-08-dis}, the proof is readily extended if the solution is  piecewise $H^2$ (in 
fact, in \cite{eym-08-dis} the situation is more complex because the Hybrid scheme is considered in the more 
general setting of the SUSHI scheme, which involves the elimination of some or all edge unknowns; see Section 
\ref{bary} below). 

\subsection{Order 2 error estimate}\label{sec-ordre2}

Under an additional assumption of existence of a specific lifting operator, compatible
with the considered Mimetic Finite Difference method, a theoretical $L^2$-error estimate of order 2
on $p$ is proved in \cite{bre-05-con}. A condition on the matrix $\mathbb{M}_E$ which ensures the
existence of such a lifting operator is given in \cite{bre-07-supercon}: in particular,
if the smallest eigenvalue of $\mathbb{U}_E$ is large enough (with respect to the inverse of
the smallest eigenvalue of $\mathbb{C}_E^T\mathbb{C}_E$ and to the largest eigenvalue
of a local inner product involving a generic lifting operator), then the
existence of a lifting operator compatible with $\mathbb{M}_E$, and thus the
super-convergence of the associated Mimetic scheme, can be proved.

Regarding the consequence on the HMMF method, this means that, if $x_K$ is
the center of gravity of $K$ and if the symmetric positive definite matrices
$(\mathbb{U}_E)_{E\in\Omega_h}$ or $(\mathbb{B}^H_K)_{K\in\mathcal M}$ 
or $(\mathbb{B}^M_K)_{K\in\mathcal M}$ are ``large enough'', then the approximate
solution $p$ given by the corresponding HMMF method converges in $L^2$ with order 2.
It does not seem easy to give a practical lower bound on
these stabilization terms which ensure that they are indeed ``large enough'' for
the theoretical proof; however, several numerical tests (both with
$x_K$ the center of gravity of $K$, or $x_K$ elsewhere inside $K$)
suggest that the HMMF methods enjoy a superconvergence property
for a wider range of parameters than those satisfying the above theoretical assumptions.


\section{Links with other methods}\label{sec-links}

We now show that, under one of its three forms and according to the choice of its parameters,
a scheme from the HMMF (Hybrid Mimetic Mixed Family) may be interpreted
as a nonconforming Finite Element scheme, as a mixed Finite Element scheme,
or as the classical  two-point flux finite volume method.

\subsection{A nonconforming Finite Element method}

In this section, we aim at identifying a hybrid FV scheme of the HMMF  with a nonconforming Finite Element method. 
Hence we use the notations of Definition \ref{def-genhybrid} and Remark \ref{rem-defhfv2}.
For any $K \in \mathcal M$ and $\sigma \in \mathcal E_K$, we  denote by $\bigtriangleup_{K,\sigma}$ the cone with vertex $x_K$ and basis $\sigma$.
For any given $\widetilde{p} \in \Xmb$ we define the piecewise linear function $\widehat p: \Omega \to \R$ by:
\[
\forall K\in\mathcal M\,,\;\forall\sigma\in\mathcal E_K\,,
\mbox{ for a.e. $x\in \bigtriangleup_{K,\sigma}$}\,:\,
\widehat p(x) = p_K + \left(\nabla_K \ploc{K} + \frac {\beta_{K}} {d_{K,\sigma}} S_{K,\sigma}(\ploc{K})\n_{K,\sigma}\right)\cdot (x - x_K)
\]
where $\beta_{K}>0$ (note that, since $\n_{K,\sigma}\cdot(\xs-x_K)=d_{K,\sigma}$,
in the particular case where $\beta_{K} = 1$ we have
$\widehat p(\xs) = p_\sigma$). We let $\widehat{H}_{\mathcal M}=\{\widehat p\,,\;\widetilde{p}\in \Xmb\}$
and we define the ``broken gradient'' of $\widehat p\in\widehat{H}_{\mathcal M}$ by

\[
\forall K\in\mathcal M\,,\;\forall
\sigma\in\mathcal E_K\,,\mbox{ for a.e. } x\in \bigtriangleup_{K,\sigma}\,:\;
\widehat \nabla  \widehat p(x) = \nabla \widehat p(x) =  \nabla_K \ploc{K} + \frac {\beta_{K}} {d_{K,\sigma}} S_{K,\sigma}(\ploc{K})\n_{K,\sigma}.
\]
We consider the following nonconforming finite element problem:
find $\widehat p\in\widehat{H}_{\mathcal M}$ such that
\begin{equation}
\forall \widehat q\in \widehat{H}_{\mathcal M}\,:\;\int_\Omega 
\widehat\Lambda(x) \widehat \nabla\widehat p(x) \cdot\widehat \nabla \widehat q(x) \dx = 
\int_\Omega f(x)  \widehat q(x) \dx,
\label{ncfe}\end{equation}
where $\widehat\Lambda(x)$ is equal to $\Lambda_K$ for a.e. $x\in K$. Then the above method
leads to a matrix which belongs to the family of matrices corresponding to the Generalized Hybrid method (and thus
also to the Generalized Mimetic and Modified Mixed methods).
Therefore, since $|\bigtriangleup_{K,\sigma}|=\frac{|\sigma|d_{K,\sigma}}{d}$, we get from \eqref{consvhf}:
\[
\int_\Omega 
\widehat\Lambda(x) \widehat \nabla\widehat p(x) \cdot\widehat \nabla \widehat q(x) \dx =
\sum_{K\in\mathcal M}|K|\Lambda_K \nabla_K \ploc{K}\cdot\nabla_K \qloc{K}
+\sum_{K\in\mathcal M} S_K(\qloc{K})^T
\mathbb{B}^H_K S_{K}(\ploc{K}),
\]
with
\[
\mathbb{B}^H_{K,\sigma,\sigma} = \frac {|\sigma|(\beta_{K})^2} {d \ d_{K,\sigma}}  \Lambda_K\n_{K,\sigma} \cdot \n_{K,\sigma},
\]
and, for $\sigma\neq \sigma'$,
\[ 
\mathbb{B}^H_{K,\sigma,\sigma'} = 0.
\]
Note  that this definition for $\mathbb{B}^H_{K}$ fulfills \refe{c1c2} under the regularity hypothesis \eqref{reghyb}, hence ensuring convergence properties which can easily be extended to this
nonconforming Finite Element method \refe{ncfe} even though it is not completely identical to a  scheme of the HMMF, since the right-hand sides do not coincide: in general,
\[
\int_\Omega f(x)  \widehat q(x) \dx \neq \sum_{K\in\mathcal M}q_K\int_K f.
\]
Indeed, this does not prevent the study of convergence since the difference between the two right-hand-sides is of order $h$.

\subsection{Relation with Mixed Finite Element methods}\label{sec-mfemethod}

In this section, we aim at identifying a particular scheme of the HMMF, under its mixed version with a Mixed Finite Element method.
Let us first recall the remark provided in \cite[Section 5.1]{bre-05-con} and
\cite[Example 1]{man-08-rt0}: on a simplicial mesh (triangular if $d=2$, tetrahedral if $d=3$),
the Raviart-Thomas RT0 Mixed Finite Element method fulfills properties \refe{S1} and \refe{S2},
and it is therefore possible to include this Mixed Finite Element scheme in the framework of the HMMF. 
Our purpose is to show that our framework also provides a Mixed Finite Element method on a general mesh
(note that, even in the case of simplices, this method differs from the Raviart-Thomas RT0 method).
We use here the notations provided by Definition \ref{def-modmixed}.
Still denoting $\bigtriangleup_{K,\sigma}$ the cone with vertex $x_K$ and basis $\sigma$,
we define, for $F\in \mathcal F$,
\begin{equation}
\forall K\in\mathcal M\,,\;\forall\sigma\in\mathcal E_K\,,\;\forall x\in \bigtriangleup_{K,\sigma}\,:\;
\widehat F_{K,\sigma}(x) := -\Lambda_K \v_K(F_K) + T_{K,\sigma}(F_K)\frac {x - x_K} {d_{K,\sigma}}.
\label{eqhat1}
\end{equation}
If $x$ belongs to the interface $\partial \bigtriangleup_{K,\sigma}\cap \partial \bigtriangleup_{K,\sigma'}$
between two cones of a same control volume,
we have $(x-x_K)\cdot\mathbf{n}_{\partial \bigtriangleup_{K,\sigma}\cap \partial \bigtriangleup_{K,\sigma'}}=0$
and thus the normal fluxes of $\widehat{F}_{K,\sigma}$
are conservative through such interfaces; moreover, for all $\sigma\in\mathcal E_K$
and all $x\in\sigma$, we have $(x-x_K)\cdot\n_{K,\sigma}=d_{K,\sigma}$ and thus,
by \refe{defpenmfv}, $\widehat F_{K,\sigma}(x)\cdot \n_{K,\sigma} = F_{K,\sigma}$;
since the elements of  $\mathcal F$ satisfy \refe{conservativite},
these observations show that the function $\widehat F$, defined by
\begin{equation}
\forall K\in\mathcal M\,,\;\forall\sigma\in\mathcal E_K\,,\mbox{ for
a.e. $x\in \bigtriangleup_{K,\sigma}$}\,:\;
\widehat F(x) = \widehat F_{K,\sigma}(x),
\label{eqhat2}\end{equation}
satisfies $\widehat F \in H_{\rm div}(\Omega)$.
Noting that
\begin{equation}
\sum_{\sigma\in\mathcal E_K} T_{K,\sigma}(F_K)\int_{\bigtriangleup_{K,\sigma}}\frac {x - x_K} {d_{K,\sigma}} \dx =
\sum_{\sigma\in\mathcal E_K} T_{K,\sigma}(F_K) \frac {|\sigma| (\xs - x_K)}{d+1} = 0
\label{sommenulle}
\end{equation}
(thanks to \refe{defvmfv}, \refe{defpenmfv} and \refe{formmag}), we have,
with $\widehat{\Lambda}$  again denoting the piecewise-constant function equal to
$\Lambda_K$ in $K$,
\begin{equation}
\int_K \widehat{\Lambda}(x)^{-1} \widehat F(x) \cdot \widehat G(x) \dx =
|K|\v_K(F_K)\cdot\Lambda_K\v_K(G_K)+\sum_{\sigma\in\mathcal E_K}
\gamma_{K,\sigma} T_{K,\sigma}(F_K)T_{K,\sigma}(G_K)
\label{goldencoin}\end{equation}
with
\[
\gamma_{K,\sigma} =\int_{\bigtriangleup_{K,\sigma}} \widehat{\Lambda}(x)^{-1}\frac {x - x_K} {d_{K,\sigma}}\cdot
\frac {x - x_K} {d_{K,\sigma}} \dx>0.
\]
The right-hand side of \refe{goldencoin} defines an inner product $\langle\cdot,\cdot\rangle_K$
which enters the framework defined by \refe{defpslocalmfv}, setting
\begin{equation}
\mathbb{B}^M_{K,\sigma,\sigma} = \gamma_{K,\sigma}, \hbox{ and }\mathbb{B}^M_{K,\sigma,\sigma'}=0
\hbox{ for }\sigma\neq \sigma',
\label{choixstab-mfe}\end{equation}
which fulfills \refe{cond-matmfv} under the regularity hypothesis \eqref{reghyb}
(hence ensuring convergence properties).
Therefore the form \listeeq{\refe{defmfd1},\refe{defmfd2}} of the resulting HMMF scheme resumes to the
following Mixed Finite Element formulation: find $(p,\widehat{F})\in \Xm \times \widehat{\mathcal F}$ such that 
\begin{eqnarray}
&&\forall  \widehat G \in \widehat{\mathcal F}\,:\;
\int_\Omega \widehat{\Lambda}(x)^{-1} \widehat F(x) \cdot \widehat G(x) \dx -
\int_\Omega p(x) \div \widehat G(x) \dx = 0\,,\label{mfe1}\\
&&\forall  q \in \Xm\,:\;
\int_\Omega q(x) \div \widehat F(x) \dx = \int_\Omega q(x) f(x) \dx\,,\label{mfe2}
\end{eqnarray}
where $\widehat{\mathcal F} = \{\widehat F\,,\;F\in {\mathcal F}\}$. Indeed,
since $|\bigtriangleup_{K,\sigma}|=\frac{|\sigma|d_{K,\sigma}}{d}$ and $\sum_{\sigma\in\mathcal E_K}|\sigma|\n_{K,\sigma}=0$,
we have, for all $F\in\widehat{\mathcal F}$,
\begin{equation}
\int_K \div \widehat{F}(x) \dx=\sum_{\sigma\in\mathcal E_K}|\sigma|T_{K,\sigma}(F_K)
=\sum_{\sigma\in\mathcal E_K}|\sigma|F_{K,\sigma}
\label{sommedivergence}\end{equation}
and \refe{mfe2} written on the canonical basis of $\Xm$ is exactly \refe{bilanmfv}.
Summing \refe{lienpFmfv} on $K$, the terms involving $p_\sigma$ vanish thanks to the conservativity of $G$
and we get \refe{mfe1}. Reciprocally, to pass from \refe{mfe1} to \refe{lienpFmfv},
one has to get back the edge values $p_\sigma$, which can be done
exactly as for the hybridization of the Generalized Mimetic method in
Section \ref{sec-hybrider_mimetic} (using $G=G(K,\sigma)$ such that
$G(K,\sigma)_{K,\sigma}=1$, $G(K,\sigma)_{L,\sigma}=-1$ if $L$
is the control volume on the other side of $\sigma$
and $G(K,\sigma)_{Z,\theta}=0$ for other control volumes $Z$ and/or edges $\theta$).

\medskip

An important element of study of the standard Mimetic method seems to be
the existence of a suitable lifting operator, re-constructing a flux unknown inside
each grid cell from the fluxes unknowns on the boundary of the grid cell (see \cite{bre-05-con}
and Section \ref{sec-ordre2}).
We claim that, for the Generalized Mimetic method corresponding to the choice \eqref{choixstab-mfe},
the flux $\widehat{F}$ given by \listeeq{\eqref{eqhat1},\eqref{eqhat2}}
provides a (nearly) suitable lifting operator~: it does not 
completely satisfy the assumptions demanded in \cite[Theorem 5.1]{bre-05-con}, but
enough so that the conclusion of this theorem still holds.

\bp For $F_K\in\mathcal F_K$, let $\widehat{F}_K$ be the restriction to $K$
of $\widehat{F}$ defined by \eqref{eqhat2}.
Then the operator $F_K\in\mathcal F_K\mapsto \widehat{F}_K\in L^2(K)$
satisfies the following properties:
\begin{equation}
\forall F_K\in\mathcal F_K\,,\;
\forall \sigma\in\mathcal E_K\,,\;\forall x\in\sigma\,:\; \widehat{F}_K(x)\cdot\n_{K,\sigma}=F_{K,\sigma}\,,
\label{lift1}\end{equation}
\begin{equation}
\forall F_K\in\mathcal F_K\,,\;\forall q\mbox{ affine function}\,:\; \int_K q(x)\div(\widehat{F}_K)(x)\dx
=\int_E q(x)\divmfd(F_K)w_E(x)\dx\,,
\label{lift2}\end{equation}
\begin{equation}
\forall F\in\R^d\,,\mbox{defining $F_K=(F\cdot\n_{K,\sigma})_{\sigma\in\mathcal E_K}$}\,:\; \widehat{F}_K=F\,,
\label{lift3}\end{equation}
for any $w_E$ satisfying \eqref{hypowE}.
\label{prop-lift}\ep

\br Properties \eqref{lift1} and \eqref{lift3} are the same as in \cite[Theorem 5.1]{bre-05-con},
but Property \eqref{lift2} is replaced in this reference by the stronger form
``$\div(\widehat{F}_K)=\divmfd(F_K)$ on $K$'' ($x_K$ is also taken as the center of gravity,
which corresponds to $w_E=1$ in \eqref{lift2}). 
\er

\bpr{ of Proposition \ref{prop-lift}}

We already noticed \eqref{lift1} (consequence of \eqref{defpenmfv} and the fact that $(x-x_K) \cdot\n_{K,\sigma}=d_{K,\sigma}$ for all $x\in\sigma$). 
If $F\in\R^d$, and $F_K=(F\cdot\n_{K,\sigma})_{\sigma\in\mathcal E_K}$, then \eqref{formmag} and
\eqref{defvmfv} show that $\Lambda_K\v_K(F_K)=-F$ and thus that $T_K(F_K)=0$,
in which case $\widehat{F}_K=-\Lambda_K\v_K(F_K)=F$ and \eqref{lift3} holds.

Let us now turn to \eqref{lift2}. For $q\equiv 1$, this relation is simply \eqref{sommedivergence}.
If $q(x)=x$, then
\[
\int_K q(x)\div(\widehat{F}_K)(x)\dx=\sum_{\sigma\in\mathcal E_K}
T_{K,\sigma}(F_K)d \int_{\bigtriangleup_{K,\sigma}}\frac{x}{d_{K,\sigma}}\dx
\]
and \eqref{sommenulle} then gives
\begin{eqnarray*}
\int_K q(x)\div(\widehat{F}_K)(x)\dx&=&\sum_{\sigma\in\mathcal E_K}
T_{K,\sigma}(F_K)d  \frac{x_K}{d_{K,\sigma}} |\bigtriangleup_{K,\sigma}|\\
&=&\sum_{\sigma\in\mathcal E_K}|\sigma| T_{K,\sigma}(F_K)x_K\\
&=&\left(\sum_{\sigma\in\mathcal E_K}|\sigma| F_{K,\sigma}\right)x_K\\
&=&\int_K q(x) \divmfd(F_K)w_E(x)\dx
\end{eqnarray*}
by assumption on $w_E$. \epr

\subsection{Two-point flux cases}\label{sec-twopoints}

We now consider isotropic diffusion tensors: 
\begin{equation}
\Lambda=\lambda(x) {\rm Id}, \label{isotrope}
\end{equation}
with $\lambda(x)\in\R$ and  exhibit cases in which the HMMF provides two-point fluxes, in the sense that the fluxes satisfy $F_E^{e} = \tau_{E,E'} (p_E - p_{E'})$, in the case where $e$ is the common edge of two neighboring cells $E$ and $E'$, with $\tau_{E,E'} \ge 0$ only depending on the grid (not on the unknowns).
Recall that using a two-point flux scheme yields a matrix with positive inverse, and  is the easiest way to ensure that, in the case of a linear scheme, monotony and local maximum principle hold.
The two-point flux scheme is  also probably the cheapest scheme in terms of implementation and computing cost. 
Moreover, although no theoretical proof is yet known, numerical evidence shows the order 2 convergence \cite{her-97-com} of the two-point scheme on triangular meshes,  taking for $x_K$ the intersection of the orthogonal bisectors \cite{her-95-err,book}.
Therefore, whenever possible, one should strive to recover this scheme when using admissible meshes (in the sense of \cite{book} previously mentioned). 
In the HMMF framework, two-point fluxes are obtained, using the notations provided by Definition \ref{sec-mfd}, if the matrix of the bilinear form $[F,G]_{X^h}$ in \refe{defmfd1} is diagonal. 
This implies, in the case of meshes such that two neighboring grid cells only have  one edge in common, that the matrices ${\mathbb M}_E$ defining the local inner product \eqref{localscal} are diagonal. 
If the  matrix ${\mathbb M}_E$ is diagonal, we get from the property ${\mathbb M}_E {\mathbb N}_E = {\mathbb R}_E$ and from \eqref{isotrope} that there exists  $\mu_E^e\in \R$ such that
\begin{equation}
\mu_E^e \n_E^e = \bar x_{e}- x_E.
\label{supadmmesh}
\end{equation}
This implies that $x_E$ is, for any face $e$ of $E$, a point of the orthogonal line to $e$ passing through $\bar x_{e}$ ($x_E$ is then necessarily unique) and that $\mu_E^e$ is the orthogonal distance between $x_E$ and $e$.
In the case of a triangle, $x_E$ is thus the intersection of the orthogonal bisectors of the sides of the triangle, which is the center of gravity only if the triangle is equilateral; hence , except in this restricted case, the original Mimetic method cannot yield a two-point flux method. 
Note also that there are meshes such that the orthogonal bisectors of the faces do not intersect, but for which nevertheless some ``centers'' in the cells exist and are such that the line joining the centers of two neighboring cells is orthogonal to their common face: these meshes are referred to ``admissible" meshes in \cite[Definition 9.1 p. 762]{book}; a classical example of such admissible meshes are the general Voronoi meshes. 
On such admissible meshes, the HMMF does not provide a two-point flux scheme for isotropic diffusion operators, although a two-point flux  Finite Volume scheme can be defined, with the desired convergence properties (see \cite{book}).

In this section, we shall call  ``super-admissible discretizations'' the discretizations which fulfill the property \refe{supadmmesh} for some choice of $(x_E)_{E\in\Omega_h}$.
We wish to show that  for all super-admissible discretizations and in the isotropic case  \eqref{isotrope}, the HMMF provides a  two-point flux scheme.
Using the notations of Definitions \ref{def-genhybrid} and  \ref{def-modmixed},  \refe{supadmmesh} is written $\n_{K,\sigma} = (\xs- x_K)/d_{K,\sigma}$ and defines our choice of parameters $(x_K)_{K\in\mathcal M}$.

\medskip

Let us take $\alpha_{K,\sigma} = \lambda_K$ in the Hybrid presentation \refe{deffluxhfv1}, denoting
by $\lambda_K$ the mean value of the function $\lambda(x)$ in $K$.
Using Definition \refe{defgradhfv} for $\nabla_K \ploc{K}$ and thanks to \refe{formmag2}, 
a simple calculation shows that
\[
\sum_{\sigma\in\mathcal E_K}\frac{|\sigma|\lambda_K} {d_{K,\sigma}}
{S}_{K,\sigma}(\ploc{K}){S}_{K,\sigma}(\qloc{K}) = 
\sum_{\sigma\in\mathcal E_K}\frac{|\sigma|\lambda_K} {d_{K,\sigma}}(p_K-p_\sigma)(q_K-q_\sigma)
- |K| \lambda_K\nabla_K\ploc{K}\cdot\nabla_K \qloc{K}.
\]
Hence we deduce from  \refe{deffluxhfv1} that $F_{K,\sigma} = \frac {\lambda_K}{d_{K,\sigma}} (p_K - p_\sigma)$
and the conservativity \refe{conservativite} leads, for an internal edge between the control volumes $K$ and $L$,
to $p_\sigma = \frac {d_{K,\sigma}\lambda_L  p_L + d_{L,\sigma} \lambda_K p_K} {d_{K,\sigma}\lambda_L  + d_{L,\sigma}\lambda_K }$.
The resulting expression for the flux becomes
\[
F_{K,\sigma} = \frac{\lambda_K\lambda_L} {d_{K,\sigma}\lambda_L +  d_{L,\sigma}\lambda_K}(p_K - p_L),
\]
which is the expression of the flux for the standard 2-points Finite Volume scheme
with harmonic averaging of the diffusion coefficient.

\medskip

It is also easy to find back this expression from the Mixed presentation \refe{defpslocalmfv},
taking
\[
\mathbb{B}^M_{K,\sigma,\sigma} = |\sigma| \frac {d_{K,\sigma}} {\lambda_K},
\]
and, for $\sigma\neq \sigma'$,
\[
\mathbb{B}^M_{K,\sigma,\sigma'} = 0.
\]
The property
\[
\sum_{\sigma\in\mathcal E_K} |\sigma|\frac {d_{K,\sigma}} {\lambda_K} T_{K,\sigma}(F_K)T_{K,\sigma}(G_K)= 
\sum_{\sigma\in\mathcal E_K} |\sigma|\frac {d_{K,\sigma}} {\lambda_K}F_{K,\sigma} G_{K,\sigma} - |K|\v_K(F_K)\cdot\lambda_K \v_K(G_K),
\]
which results from \refe{defvmfv}, \refe{defpenmfv}, \refe{formmag2} and \refe{supadmmesh}
(under the form $d_{K,\sigma}\n_{K,\sigma}=\xs-x_K$), shows that
\[
\langle F_K,G_K\rangle_K = \sum_{\sigma\in\mathcal E_K} |\sigma|\frac{d_{K,\sigma}}{\lambda_K}F_{K,\sigma} G_{K,\sigma}.
\]
Thanks to  \refe{lienpFmfv}, this gives $F_{K,\sigma} = \frac {\lambda_K}{d_{K,\sigma}} (p_K - p_\sigma)$
and we conclude as above.

\subsection{Elimination of some edge unknowns}\label{bary}

In the study of the hybrid version  \refe{defhfv2} of  the HMMF \cite{eym-08-ben,eym-08-dis}, it was suggested to replace the space $\Xmbzero$ by the space $\XmbzeroB$ defined by $\XmbzeroB = \{ \widetilde{p}\in \Xmbzero$ such that $p_\sigma = \sum_{K\in{\mathcal M}_\sigma} \beta_\sigma^K p_K$ if $\sigma\in{\mathcal B}\}$, where:
\begin{itemize}
\item[i)] ${\mathcal M}_\sigma$ is a subset of ${\mathcal M}$, including a few cells 
(in general, less than $d+1$, where we recall that $d$ is the dimension of the space)
``close'' from $\sigma$,
\item[ii)] the coefficients $\beta_\sigma^K$ are barycentric weights of the point $\bar x_\sigma$
with respect to the points $x_K$, which means that $\sum_{K\in{\mathcal M}_\sigma} \beta_\sigma^K = 1$ and
$\sum_{K\in{\mathcal M}_\sigma} \beta_\sigma^K x_K= \bar x_\sigma$,
\item[iii)] ${\mathcal B}$ is any subset of the set of all internal edges (the cases of the empty set
or of the full set itself being not excluded).
\end{itemize}
Then the scheme is defined by: find $\widetilde{p}\in \XmbzeroB$ such that
\begin{equation}
\forall \widetilde{q}\in \XmbzeroB\,:\;
\sum_{K\in\mathcal M}|K|\Lambda_K \nabla_K \ploc{K}\cdot\nabla_K \qloc{K}
+\sum_{K\in\mathcal M} S_K(\qloc{K})^T
\mathbb{B}^H_K S_{K}(\ploc{K})
=\sum_{K\in\mathcal M}q_K\int_K f.
\label{defhfv2bary}\end{equation}
In the case where ${\mathcal B} = \emptyset$, then the method belongs to the HMMF and in the case where ${\mathcal B}$ is the full set of the internal edges, then there is no more
edge unknowns, and we get back a cell-centered scheme. In the intermediate cases, we get schemes
where the unknowns are all the cell unknowns, and the edge unknowns $p_\sigma$ with
$\sigma\notin {\mathcal B}$.

This technique has been shown in \cite{eym-08-ben} and \cite{eym-08-dis} to fulfill the convergence and error estimates requirements in the case of diagonal matrices $\mathbb{B}^H_K$. 
It can be applied, with the same convergence properties, to the HMMF with symmetric positive definite matrices $\mathbb{B}^H_K$ as in Section \ref{sec-convergence}.

\section{Appendix}

\subsection{About the Generalized Mimetic definition}\label{appen-genmim}

We prove here two results linked with the definition of the Generalized Mimetic
method: the existence of a weight function satisfying \refe{hypowE}
and the equivalence between \refe{GS2} and \refe{choixmatmfd}.

\bl If $E$ is a bounded non-empty open subset of $\R^d$ and $x_E\in\R^d$, then there exists
an affine function $w_E:\R^d\to \R$ satisfying \refe{hypowE}.
\label{exist-poids}\el

\bpr{ of Lemma \ref{exist-poids}}

We look for $\xi\in\R^d$ such that $w_E(x)=1+\xi\cdot(x-\bar x_E)=1+(x-\bar x_E)^T\xi$
satisfies the properties (where $\bar x_E$ is the center of gravity of $E$).
The first property of \refe{hypowE} is straightforward since $\int_E(x-x_E)\dx=0$, and the
second property is equivalent to $|E|\bar x_E+\int_E x(x-\bar x_E)^T\xi\dx=|E|x_E$;
since $\int_E \bar x_E(x-\bar x_E)^T\dx=0$, this boils down to
\begin{equation}
\left(\int_E (x-\bar x_E)(x-\bar x_E)^T\dx\right)\xi = |E|(x_E-\bar x_E).
\label{rel-xi}\end{equation}
Let $J_E$ be the $d\times d$ matrix $\int_E (x-\bar x_E)(x-\bar x_E)^T\dx$: we have, for all 
$\eta\in\R^d\backslash\{0\}$,
$J_E\eta\cdot \eta = \int_E ((x-\bar x_E)\cdot\eta)^2\dx$ and the function
$x\to (x-\bar x_E)\cdot\eta$ vanishes only on an hyperplane of $\R^d$; this proves
that $J_E\eta\cdot\eta>0$ for all $\eta\not=0$. Hence, $J_E$ is
invertible and there exists  (a unique)  $\xi$ satisfying \refe{rel-xi}, which
concludes the proof. \epr

\bl Let $[\cdot,\cdot]_E$ be a local inner product on the space of the fluxes
unknowns of a grid cell $E$, and let $\mathbb{M}_E$ be its matrix. Then $[\cdot,\cdot]_E$
satisfies \eqref{GS2} (with $w_E$ satisfying \refe{hypowE}) if and only if
$\mathbb{M}_E$ satisfies \refe{choixmatmfd} (with $\mathbb{R}_E$ defined by
\refe{defnewR} and $(\mathbb{C}_E,\mathbb{U}_E)$ defined by \eqref{mf2}  and  \eqref{mf3}).
\label{lem-equivGS2}\el

\bpr{ of Lemma \ref{lem-equivGS2}}

It is known \cite{bre-05-fam} that, for the standard Mimetic method, \eqref{S2} is
equivalent to $\mathbb{M}_E\mathbb{N}_E=\bar{\mathbb{R}}_E$ with $\bar{\mathbb{R}}_E$ defined by
\eqref{mf1} and $\mathbb{N}_E$ defined in \eqref{mf2};
similarly, it is quite easy to see that \eqref{GS2} (with $w_E$ satisfying
\eqref{hypowE}) is equivalent to
\begin{equation}
\mathbb{M}_E\mathbb{N}_E=\mathbb{R}_E
\label{GS2-mat}\end{equation}
with $\mathbb{R}_E$ defined by \refe{defnewR}: indeed, \eqref{GS2} with $q=1$ is simply the
definition \refe{defdivmfd} of the discrete divergence operator (because $\int_E w_E(x)\dx=|E|$) and, with $q(x)=x_j$, since $\int_E x_jw_E(x)\dx=|E|(x_E)_j$, \eqref{GS2} boils down to
\[
G_E^T \mathbb{M}_E(\mathbb{N}_E)_j+\sum_{i=1}^{k_E}G_E^{e_i}|e_i|(x_E)_j=
\sum_{i=1}^{k_E}G^{e_i}_E|e_i|(\bar x_{e_i})_j,
\]
which is precisely $G_E^T \mathbb{M}_E(\mathbb{N}_E)_j=G_E^T(\mathbb{R}_E)_j$
with $(\mathbb{R}_E)_j$ the $j$-th column of $\mathbb{R}_E$. We therefore only need to compare
\eqref{choixmatmfd} with \eqref{GS2-mat}.

\medskip

Let us first assume that $\mathbb{M}_E$ satisfies \eqref{choixmatmfd}. The generic formula \eqref{formmag}
implies
\begin{equation}
\mathbb{R}_E^T(\mathbb{N}_E)_j=|E|(\Lambda_E)_j
\label{RN}\end{equation}
and thus
\begin{equation}
\frac{1}{|E|}\mathbb{R}_E\Lambda_E^{-1}\mathbb{R}_E^T\mathbb{N}_E=\mathbb{R}_E.
\label{RLambdaN}\end{equation}
Since $\mathbb{C}_E^T\mathbb{N}_E=0$ by definition of $\mathbb{C}_E$,
this shows that $\mathbb{M}_E$ satisfies \eqref{GS2-mat}.

\medskip

Let us now assume that $\mathbb{M}_E$ satisfies \eqref{GS2-mat} and let us consider
the symmetric matrix
\begin{equation}
\widetilde{\mathbb{M}}_E=\mathbb{M}_E-\frac{1}{|E|}\mathbb{R}_E\Lambda_E^{-1}\mathbb{R}_E^T.
\label{Mtilde}\end{equation}
By \eqref{GS2-mat} and \eqref{RLambdaN}, we have $\widetilde{\mathbb{M}}_E\mathbb{N}_E=0$,
and the columns of $\mathbb{N}_E$ are therefore in the kernel of $\widetilde{\mathbb{M}}_E$.
The definition of $\mathbb{C}_E$ shows that the columns of $\mathbb{N}_E$ span the
kernel of $\mathbb{C}_E^T$: we therefore have $\ker(\mathbb{C}_E^T)\subset \ker(\widetilde{\mathbb{M}}_E)$
and we deduce the existence of a $k_E\times (k_E-d)$ matrix $\mathbb{A}$ such that 
\begin{equation}
\widetilde{\mathbb{M}}_E=\mathbb{A}\mathbb{C}_E^T.
\label{Mtilde2}\end{equation}
By symmetry of $\widetilde{\mathbb{M}}_E$ we have $\mathbb{A}\mathbb{C}_E^T=\mathbb{C}_E\mathbb{A}^T$
and thus $\mathbb{A}\mathbb{C}_E^T\mathbb{C}_E=\mathbb{C}_E\mathbb{A}^T\mathbb{C}_E$; but
$\mathbb{C}_E^T\mathbb{C}_E$ is an invertible $(k_E-d)\times (k_E-d)$ matrix (since $\mathbb{C}_E$
is of rank $k_E-d$) and therefore $\mathbb{A}=\mathbb{C}_E\mathbb{A}^T\mathbb{C}_E
(\mathbb{C}_E^T\mathbb{C}_E)^{-1}=\mathbb{C}_E\mathbb{U}_E$ for some $(k_E-d)\times (k_E-d)$
matrix $\mathbb{U}_E$. Gathering this result with \eqref{Mtilde} and \eqref{Mtilde2}, we
have proved that $\mathbb{M}_E$ satisfies \eqref{choixmatmfd} for some
$\mathbb{U}_E$, and it remains to prove that this last matrix is symmetric definite positive
to conclude. By \eqref{choixmatmfd} and the symmetry of $\mathbb{M}_E$ and
$\frac{1}{|E|}\mathbb{R}_E\Lambda_E^{-1}\mathbb{R}_E^T$
we have $\mathbb{C}_E\mathbb{U}_E\mathbb{C}_E^T=\mathbb{C}_E\mathbb{U}_E^T\mathbb{C}_E^T$
and, since $\mathbb{C}_E^T$ is onto and $\mathbb{C}_E$ is one-to-one, we deduce
$\mathbb{U}_E=\mathbb{U}_E^T$. To prove that $\mathbb{U}_E$ is definite positive,
we use \eqref{choixmatmfd} and the fact that $\mathbb{M}_E$ is definite positive
to write
\begin{equation}
\forall\xi\in\ker(\mathbb{R}_E^T)\,,\;\xi\not=0\,:\;(\mathbb{U}_E\mathbb{C}_E^T\xi)\cdot (\mathbb{C}_E^T\xi)>0.
\label{Usdp}\end{equation}
This shows in particular that $\ker(\mathbb{R}_E^T)\cap \ker(\mathbb{C}_E^T)=\{0\}$ and thus,
since $\ker(\mathbb{R}_E^T)$ has dimension $k_E-d$, that the image by
$\mathbb{C}_E^T$ of $\ker(\mathbb{R}_E^T)$ is $\R^{k_E-d}$. Equation \eqref{Usdp} then proves
that $\mathbb{U}_E$ is definite positive on the whole of $\R^{k_E-d}$ and
the proof is complete. \epr

\subsection{An algebraic lemma}

\bl
Let $X$, $Y$ and $Z$ be finite dimension vector spaces and $A:X\to Y$,
$B:X\to Z$ be two linear mappings with identical kernel.
Then, for all inner product $\{\cdot,\cdot\}_Y$ on $Y$, there exists
an inner product $\{\cdot,\cdot\}_Z$ on $Z$ such that,
for all $(x,x')\in X^2$, $\{Bx,Bx'\}_Z=\{Ax,Ax'\}_Y$.
\label{lem-clef}\el

\bpr{ of Lemma \ref{lem-clef}}

Let $N=\ker(A)=\ker(B)$. The mappings $A$ and $B$ define
one-to-one mappings $\bar A:X/N\to Y$
and $\bar B:X/N\to Z$ such that, if $\bar x$ is the class of $x$,
$Ax=\bar A\bar x$ and $Bx=\bar B\bar x$. We can therefore work
with $\bar A$ and $\bar B$ on $X/N$ rather than with $A$ and $B$ on
$X$, and assume in fact that $A$ and $B$ are one-to-one.

Then $A:X\to {\rm Im}(A)$ and $B:X\to {\rm Im}(B)$ are isomorphisms
and, if $\{\cdot,\cdot\}_Y$ is an inner product on $Y$,
we can define the inner product $\{\cdot,\cdot\}_{{\rm Im}(B)}$ on ${\rm Im}(B)$ the following
way: for all $z,z'\in {\rm Im}(B)$, $\{z,z'\}_{{\rm Im}(B)}=\{AB^{-1}z,AB^{-1}z'\}_Y$
(this means that $\{Bx,Bx'\}_{{\rm Im}(B)}=\{Ax,Ax'\}_Y$ for all $x,x'\in X$). This inner product
is only defined on ${\rm Im}(B)$, but we extend it
to $Z$ by choosing $W$ such that ${\rm Im}(B)\oplus W=Z$,
by taking any inner product $\{\cdot,\cdot\}_W$ on $W$ and by
letting $\{z,z'\}_Z=\{z_B,z'_B\}_{{\rm Im}(B)}+\{z_W,z'_W\}_W$ for all
$z=z_B+z'_W\in Z={\rm Im}(B)\oplus W$ and $z'=z'_B+z'_W\in Z$. This extension
of $\{\cdot,\cdot\}_{{\rm Im}(B)}$ preserves the property $\{Bx,Bx'\}_Z=\{Ax,Ax'\}_Y$. \epr

\br The proof gives a way to explicitly compute $\{\cdot,\cdot\}_Z$ from
$\{\cdot,\cdot\}_Y$, $A$ and $B$: find a supplemental space $G$ of $\ker(A)=\ker(B)$ in $X$,
compute an inverse of $B$ between $G$ and ${\rm Im}(B)$, deduce $\{\cdot,\cdot\}_{{\rm Im}(B)}$
and extend it to $Z$ by finding a supplemental space of ${\rm Im}(B)$.
\label{rem-explicite0}\er

\end{document}